\documentclass[10pt]{article}
\usepackage{amssymb,amsmath}
\usepackage{graphicx}
\usepackage{enumerate}
\usepackage{color}
\usepackage{multirow}
\def\vec#1{{\boldsymbol #1}}
\def\disc#1{{\mathbf #1}}
\def\calbf#1{{\boldsymbol{\mathcal #1}}}
\newcommand{\R}{{\mathbb R}}
\newcommand{\Z}{{\mathbb Z}}

\def\ds{\displaystyle} 
\newcommand{\der}{{\mathrm{d}}}
\def\step{\Delta t}

\newcommand{\x}{{\boldsymbol{x}}}
\newcommand{\I}{{\mathrm{Id}}}
\newcommand{\Ibf}{{\mathbf{Id}}}
\newcommand{\ie}{{{\it ,\,i.e.,\,}}}
\newcommand{\eg}{{{\it ,\,e.g.,\,}}}
\newcommand{\MR}{{\mathcal M}}
\newcommand{\thr}{{\mathcal T}}
\newcommand{\adap}{{\mathcal A}}
\newcommand{\leaf}{{\mathrm{L}}}

\newcommand{\One}{{\mathbf{1}}}
\newcommand{\tolnew}{{\eta_{\mathrm{Newt}}}}
\newcommand{\knew}{{k_{\mathrm{Newt},\max}}}
\newcommand{\klinjac}{{k_{\mathrm{LS},\disc J}}}
\newcommand{\tollin}{{\eta_{\mathrm{LS}}}}
\newcommand{\toltime}{{\eta_{\mathrm{RK}}}}
\newcommand{\tolMR}{{\eta_{\mathrm{MR}}}}
\newcommand{\new}{{\mathrm{new}}}
\newcommand{\klin}{{k_{\mathrm{LS}}}}
\newcommand{\Or}{{\mathcal O}}
\newcommand{\F}{{\mathrm{F}}}
\newcommand{\Ox}{{\mathrm{O}}}
\newcommand{\e}{{\mathrm e}}
\def\tlvs{\vrule height 1em  width 0pt} 

\title{High order implicit time integration schemes
on multiresolution adaptive grids for stiff PDE's}

\author{
Max Duarte\footnotemark[1]\ \footnotemark[2]
\and
Richard R. Dobbins\footnotemark[1]
\and
Mitchell D. Smooke\footnotemark[1]
}

\begin{document}

\maketitle

\renewcommand{\thefootnote}{\fnsymbol{footnote}}

\footnotetext[1]{
Department of Mechanical Engineering, 
Yale University, Becton Center, 15 Prospect Street, New Haven, CT 06520, USA
({\tt \{richard.dobbins,mitchell.smooke\}@yale.edu}).
}

\footnotetext[2]{
{\it Present address:} 
CD-adapco, 200 Shepherds Bush Road, London W6 7NL, UK
({\tt max.duarte@cd-adapco.com}).
}

\renewcommand{\thefootnote}{\arabic{footnote}}

\begin{abstract}
We consider high order, implicit Runge--Kutta schemes
to solve time--dependent stiff PDEs
on dynamically adapted grids generated by multiresolution
analysis 
for unsteady problems disclosing localized fronts.
The multiresolution finite volume scheme yields 
highly compressed representations within a user--defined accuracy tolerance,
hence strong reductions of computational requirements to solve
large, coupled nonlinear systems of equations.
SDIRK and RadauIIA Runge--Kutta schemes are implemented with 
particular interest in those with $L-$stability properties and
accuracy--based time--stepping capabilities.
Numerical evidence is provided of the computational efficiency
of the numerical strategy to cope with highly unsteady problems
modeling various physical scenarios with a broad
spectrum of time and space scales.
\end{abstract}



\pagestyle{myheadings}
\thispagestyle{plain}
\markboth{DUARTE, DOBBINS, SMOOKE}
{TIME IMPLICIT SCHEMES ON MULTIRESOLUTION GRIDS}

\section{Introduction}\label{sec:intro}

In many scientific applications,
from biomedical models
to combustion or air pollution modeling,
stiff differential equations must be solved 
to carry out numerical simulations.
A straightforward notion of stiffness was given by Hairer \& Wanner 
by stating that {\it stiff equations are problems
for which explicit methods don't work}
\cite{Hairer96}.
The latter is mainly due to the broad spectrum of 
physical or numerical time scales that a numerical solver 
must deal with whenever a stiff time--dependent differential
equation is being solved.
Robust and stable methods are thus required for stiff problems
in order to properly handle and damp out fast transients.
In the past decades high order
implicit Runge--Kutta schemes 
with excellent stability properties
were developed and widely
investigated to solve stiff problems modeled by ODEs
(see \cite{Hairer96} \S~IV and references therein).
The same methods can be naturally considered to solve
stiff problems originating from
time--dependent PDEs discretized in space.
However, the high performance of implicit Runge--Kutta methods
for stiff ODEs
is often adversely affected by the
size of the systems of nonlinear equations
arising in the case of semi--discrete PDEs.
In particular
phenomena involving localized fronts,
as considered in this work, 
commonly require fine spatial
representations, hence potentially large
systems of equations.
Significant effort is required to
achieve numerical implementations that
solve the corresponding algebraic problems
at reasonable computational expenses in terms
of both CPU time and memory.

Low order implicit schemes were already successfully used to 
simulate very complex problems modeled by stiff PDEs.
This is the case, for instance, for the numerical simulation
of combustion flames accounting for detailed chemical kinetics
and multi--species transport
(see\eg \cite{Bennett2009,Smooke2013} 
and references therein).
The computational performance of high order Runge--Kutta methods,
implemented in well--established production aerodynamics codes,
was also assessed in the context of laminar and turbulent 
compressible flows \cite{Bijl2002,Carpenter05}.
In particular Jacobian--free Newton--Krylov methods were investigated
in conjunction with high order implicit schemes \cite{Zingg04,Bijl2005}
to further reduce the computational requirements
(see\eg  \cite{Knoll2004} for a review on this subject).
Similarly,
high order space discretization
schemes have also been implemented in this context,
reducing in practice the computational stencils, 
hence the size of the nonlinear systems
\cite{Noskov2005,Noskov2007}. 
Easing the computational load is also achievable
by designing efficient parallelization techniques 
as developed, for instance,
in \cite{Tosatto2011} for reactive flow solvers.
Taking into account that
grid adaptation techniques for unsteady problems
disclosing localized fronts 
are specifically designed to yield high data compression,
we exploit this capability here to efficiently implement
implicit integration schemes for stiff PDEs.
This strategy was already adopted, for instance,
in \cite{Bennett1998,Bennett1999}, together with low order implicit solvers.

Among the many adaptive meshing approaches developed in the literature,
we consider in this work adaptive multiresolution
schemes based on \cite{Harten94,Harten95},
namely the multiresolution finite volume
scheme introduced in \cite{Cohen03}
for conservation laws.
Besides the inherent advantages
of grid adaptation, multiresolution
techniques rely on biorthogonal wavelet decomposition
and thus offer a rigorous 
mathematical framework for adaptive meshing
schemes \cite{cohen2000a,muller2003}.
Consequently,
not only approximation errors 
coming from grid adaptation
can be tracked, but general and robust 
solvers can be implemented since 
the wavelet decomposition is
independent of any physical particularity 
of the problem and accounts
only for the spatial regularity of the 
discrete variables at a given simulation time.
Adaptive multiresolution schemes have been successfully
implemented for the simulation of compressible fluids
(see\eg \cite{Brix2011,Domingues2011} 
and references therein),
as well as for the numerical solution of time--dependent,
parabolic \cite{Roussel03,Burger08} and stiff parabolic 
\cite{Duarte11_SISC,Dumont2013}
PDEs.
Nevertheless, to the best of our knowledge this is the fist attempt
to implement high order implicit time integration schemes 
in the context of the adaptive multiresolution 
finite volume method to solve stiff PDEs.

The paper is organized as follows.
We give in Section~\ref{sec:scheme}
a short introduction on multiresolution finite volume schemes
and implicit Runge--Kutta schemes,
in particular of
SDIRK-- 
(Singly Diagonally Implicit Runge--Kutta)
and RadauIIA-- type.
Some key aspects of the numerical implementation
of time implicit schemes on multiresolution grids
are detailed in Section~\ref{sec:imple}.
Finally, the numerical solution of several 
stiff time--dependent PDEs is investigated in 
Section~\ref{sec:num_res}.

\section{Numerical methodology}\label{sec:scheme}
Let us consider a parabolic, time--dependent PDE,
\begin{equation}\label{eq:gen_prob}
\left.
\begin{array}{ll}
\partial_t \vec u = \vec F(\partial^2_{\x} \vec u, \partial_{\x} \vec u, \vec u), &
t > t_0,\, \x \in \R^d,\\
\vec u(t_0,\x) = \vec u_0(\x),&
t = t_0,\, \x \in \R^d,
\end{array}
\right\}
\end{equation}
where $\vec u:\R\times \R^d \to \R^m$ and
$\vec F:\R^m \to \R^m$, for a model with $m$ variables.
For many physically inspired systems, the right hand side $\vec F$ can 
in general be written as
\begin{equation}\label{eq:def_F}
\vec F(\partial^2_{\x} \vec u, \partial_{\x} \vec u, \vec u) =
\vec F(\vec u) = \vec F_1(\vec u) + \vec F_2(\vec u) + \ldots,
\end{equation}
where the $\vec F_i(\vec u) $, $i=1,\ldots$, stand for different
physical processes.
For instance, a scalar nonlinear reaction--diffusion equation
with $u:\R\times \R^d \to \R$
would be given by $F_1(u)=- \partial_{\x} \cdot (D(u) \partial_{\x} u)$
and $F_2(u)= f(u)$ for some diffusion coefficient,
$D:\R \to \R$, and a nonlinear function, $f:\R \to \R$.

\subsection{Multiresolution analysis}\label{sec:MR}
Without loss of generality
we perform a 
finite volume discretization of problem (\ref{eq:gen_prob})
with $m=1$ for the sake of simplicity.
According to the multiresolution finite volume scheme
\cite{Cohen03},
we consider a set of nested dyadic grids 
over a computational domain
$\Omega \subset \mathbb{R}^{d}$ as follows.
Each cell $\Omega _{\lambda}$, $\lambda \in S_j$, 
is the union of $2^d$ finer cells of equal size $\Omega_{\mu}$, 
$\mu \in S_{j+1}$.
The sets $S_j$ and $S_{j+1}$ are thus consecutive embedded grids over $\Omega$,
where $j=0,1,\ldots,J$, corresponds to the grid--level, from
the coarsest to the finest grid\ie
$j$ equal to 0 and $J$, respectively.
We denote
$\disc U_j:=(u_{\lambda})_{\lambda \in S_j}$ as the 
spatial representation of $u$
on the grid $S_j$, where $u_{\lambda}$
represents the cell--average of 
$u :\, \R \times \R^d \to \R$
in $\Omega _{\lambda}$:
\begin{equation}\label{eq3:average_finite_vol}
  u_{\lambda} := |\Omega_{\lambda}|^{-1} \int_{\Omega_{\lambda}}  u(t,\x)\, \der \x, \quad \x\in \R^d.
 \end{equation}
 
Data at different levels of discretization are related
by two inter--level transformations:
the {\it projection} and {\it prediction}
operators, briefly defined in Appendix~\ref{app:multiresolution}.
Based on these operations, 
a multiresolution analysis allows one to define 
a one--to--one correspondence between two consecutive grid--levels:
\begin{equation}\label{eq3:one_one_cor}
\disc U_j\longleftrightarrow (\disc U_{j-1},\disc D_j),
\end{equation}
where the $\disc D_j$ array gathers the so--called {\it details}.
The latter can be seen as estimators of the local spatial regularity
of a given discretized function, in this case $\disc U_j$,
and represent the information lost when coarsening the spatial grid,
in this case from $\disc U_j$ to $\disc U_{j-1}$.
By iteration of this decomposition, we get a 
multi--scale representation of $\disc U_J$
in terms of $\disc M_J := (\disc U_0,\disc D_1,\disc D_2,\cdots,\disc D_J)$:
\begin{equation}\label{eq3:M_U_M}
\MR:\disc U_J\longmapsto \disc M_J,
\end{equation}
and similarly, its inverse $\MR^{-1}$.
This multi--scale transform amounts to 
a representation of $\disc U_J$
in a wavelet space spanned by a biorthogonal wavelet basis.
Further details can be found in \cite{cohen2000a,muller2003}.

While the transformation (\ref{eq3:M_U_M}) is exact and the multi--scale 
representation can be performed back and forth, real computational
benefit is achieved by introducing a {\it thresholding}
operator, as shown in Appendix~\ref{app:multiresolution}.
This operator basically discards cells of smooth regularity
whose values can be recomputed within an accuracy tolerance 
whenever needed.
As a result a multiresolution approximation $\disc U_J^{\epsilon}$
is obtained.
Defining the following normalized $\ell^2$--norm:
\begin{equation*}\label{eq3:normalized_l2}
\|\disc U_J\|_2^2:= 2^{-dJ} \ds \sum_{\lambda \in S_J} (u_{\lambda})^2,
\end{equation*}
which corresponds to the $L^2$--norm of a piecewise constant function,
it can be shown that \cite{Duarte14_Poisson}
\begin{equation}\label{eq3:adap_error_eps}
 \|\disc U_{J} - \disc U_J^{\epsilon} \|_2 \leq C \tolMR,
\end{equation}
where $\tolMR$ corresponds to an accuracy tolerance\footnote{
Bound (\ref{eq3:adap_error_eps}) was similarly shown in \cite{Cohen03} 
for both uniform and $\ell^1$--norm.}.

So much is true for steady problems.
When solving time--dependent problems,
the same behavior is expected
in terms of numerical errors introduced by the multiresolution
approximation.
(The spatially adapted grid is fixed during a given time integration step.)
The latter was mathematically proved for hyperbolic problems
in an $L^1$--norm for both 
classical and inhomogeneous 
conservation laws in
\cite{Cohen03} and \cite{Hovhannisyan2010}, respectively.
Moreover,
numerical evidence proves similar behaviors for 
time--dependent, parabolic
\cite{Roussel03,Burger08,Bendahmane09}
and stiff parabolic
\cite{Duarte11_SISC,Duarte11_JCP,DuarteCFlame} PDEs.

\subsection{Implicit Runge--Kutta schemes}\label{sec:IRK}
Let us now consider problem (\ref{eq:gen_prob})
discretized on an adapted grid obtained by multiresolution analysis:
\begin{equation}\label{eq:gen_disc_prob}
\left.
\begin{array}{ll}
\der_t \disc U = \disc F (\disc U), &
t > t_0,\\
\disc U(t_0) = \disc U_0,&
t = t_0.
\end{array}
\right\}
\end{equation}
For the ease of reading we denote
$\disc U_J^{\epsilon}$ simply as $\disc U$
of size $m\times N$, where $N$ corresponds to the number of cells
in the adapted grid
and thus $\disc F:\R^{m\times N} \to \R^{m\times N}$.

Given a time step $\step$
we consider an  implicit Runge--Kutta (IRK) scheme
of order $p$
for the numerical integration of
the semi--discretized  problem (\ref{eq:gen_disc_prob}).
An $s$--stage Runge--Kutta scheme is in general defined through a set of 
arrays $\vec b$, $\vec c \in \R^s$, such that
$\vec b=(b_1, \ldots, b_s)^T$ and $\vec c=(c_1,\ldots, c_s)^T$, and 
a matrix $\vec A \in \mathcal{M}_s(\R)$ such that 
$\vec A=(a_{ij})_{1 \leq i,j \leq s}$.
These coefficients
define the stability properties and the order
conditions of the method,
and are usually arranged in a Butcher tableau according to
\begin{equation*}
\begin{tabular}{c|c }
$\vec c$ & $\vec A$ \\
\hline \\[-2.5ex]
& $\vec b^T$
\end{tabular}.
\end{equation*}
In practice, given a set of arrays 
$\disc z_1,\ldots,\disc z_s \in \R^{m\times N}$,
we have to solve the nonlinear system 
\begin{equation}\label{eq2:nonlinear_sys}
 \left(
\begin{array}{c}
\disc z_1\\
\vdots\\
\disc z_s
\end{array}\right)
=\disc A\left(
\begin{array}{c}
\step \disc F (t_0+c_1\step,\disc U_0+ \disc z_1)\\
\vdots\\
\step \disc F (t_0+c_s\step, \disc U_0+ \disc z_s)
\end{array}\right),
\end{equation}
where $\disc A$ is a square block--matrix 
of size $s \times m \times N$
built with the coefficients $(a_{ij})_{1 \leq i,j \leq s}$
(see more details in Appendix~\ref{app:IRK}).
The solution $\disc U(t_0+\step)$
is then approximated by $\disc U_1$, computed as
\begin{equation}\label{eq2:u1_radau}
\disc U_1= \disc U_0+\ds \sum_{i=1}^sd_i \disc z_i,
\qquad
\vec d^T := (d_1,\ldots,d_s)=(b_1,\ldots,b_s)\vec A^{-1}.
\end{equation}
If all the elements of the matrix of coefficients $\vec A$
are non--zero, then we say we are considering a
{\it fully IRK scheme} \cite{Hairer96}.
Moreover, if 
\begin{equation}\label{eq:stiffly_acc}
a_{sj} = b_j, \quad j = 1,\ldots,s
\end{equation}
then the last stage corresponds to the solution $\disc U_1$ 
and thus $\vec d^T = (0,0,\ldots,0,1)$ in (\ref{eq2:u1_radau}).
Methods satisfying (\ref{eq:stiffly_acc})
are called {\it stiffly accurate} \cite{Prothero74}
and are particularly appropriate for the solution of 
(stiff) singular perturbation problems and for 
differential--algebraic equations \cite{Hairer88,Hairer96}.

An IRK approximation amounts then to 
solving a nonlinear system of equations of size $s \times m \times N$.
The latter can be achieved by considering a {\it simplified} Newton method for
system (\ref{eq2:nonlinear_sys}), recast as
\begin{equation}
\calbf G(\disc Z) := \disc Z - \step \disc A \calbf F(\disc Z)= \disc 0,
\end{equation}
where $\disc Z:=(\disc z_1,\ldots,\disc z_s)^T$ and 
$\calbf F(\disc Z):=(\disc F(t_0+c_1\step,\disc U_0+\disc z_1),\ldots,\disc F(t_0+c_s\step,\disc U_0+\disc z_s))^T$.
The $(k+1)$--th approximation of the solution $\disc Z$ is thus computed 
in two steps.
First, we solve the following linear system for $\delta \disc Z^k \in \R^{s\times m \times N}$:
\begin{equation}\label{eq2:Newton_simp}
(\Ibf_{s\times m \times N} - \step \disc A \disc J) \delta \disc Z^k 
= -\disc Z^k+ \step \disc A \calbf F(\disc Z^k),
\end{equation}
where $\disc J$ is a 
a square block--matrix 
of size $s \times m \times N$
consisting of $s$ rows and $s$ columns of size 
$m \times N$ given by the Jacobians
$\disc J_0 := \partial_{\disc U} \disc F (t_0,\disc U_0)$.
Then, the previous solution $\disc Z^k$ is corrected according to 
\begin{equation}\label{eq:fin_Newton}
\disc Z^{k+1}=\disc Z^k+\delta \disc Z^k.
\end{equation}
A standard way of initializing the iterative algorithm
considers simply 
\begin{equation}\label{eq:ini_Newton}
\disc z_i^0 =\disc U_0, \quad
i=1,\ldots,s.
\end{equation}

\subsubsection{SDIRK schemes}
An alternative to solving a large 
nonlinear system of size $s \times m \times N$
is to consider an SDIRK 
scheme where
$a_{ij}=0$ for $j>i$,
and with equal diagonal coefficients\ie
$a_{ii}=\gamma$, $i=1,\ldots,s$.
In general, $A$-- or $L$--stable SDIRK schemes
can be built of order $p \leq s+1$ or 
$p \leq s$, respectively (see \cite{Hairer96} \S~IV.6 for more details).
However, the stage order $q$ of these schemes, that is, the order achieved
by one single stage, is limited to $1$.

The main idea is thus to successively
solve the $s$ stages by considering an $m\times N$--dimensional
system at each stage, that is, for $i=1,\ldots,s$,
\begin{equation}
 \disc z_i = \step \gamma \disc F(t_0+c_i\step,\disc U_0+ \disc z_i) 
+ \step \ds \sum_{j=1}^{i-1} a_{ij}\disc F(t_0+c_j\step,\disc U_0+ \disc z_j),
\end{equation}
where the second term in the right--hand side is already known
at the current stage.
Adopting the same simplified Newton technique,
this time stage--wise,
we have to solve the following linear system
for  $\delta \disc z_i^k \in \R^{m \times N}$:
\begin{align}\label{eq2:Newton_simp_SDIRK}
(\Ibf_{m \times N} - \step \gamma \disc J_0) \delta \disc z_i^k 
= & -\disc z_i^k+ \step \gamma \disc F(t_0+c_i\step,\disc U_0+ \disc z^k_i ) \nonumber \\
& + \step \sum_{j=1}^{i-1} a_{ij}\disc F(t_0+c_j\step,\disc U_0+ \disc z_j).
\end{align}
A simple way of initializing the iterative algorithm
at each stage considers
\begin{equation}\label{eq:ini_SDIRK4}
\disc z_1^0 =\disc U_0, \quad
\disc z_i^0 =\disc U_0+ \disc z_{i-1}, \quad
i=2,\ldots,s.
\end{equation}
We will denote as SDIRK2, SDIRK3, and SDIRK4,
respectively, the second, third, and fourth order SDIRK schemes
here considered (Butcher tableaux in Appendix~\ref{app:Butcher_tab}).

\subsubsection{RadauIIA schemes}
Fully IRK schemes
with a number of stages below its approximation order
can be built
based on collocation methods
\cite{Guillon69,Wright71},
together with the simplified order conditions
introduced by Butcher \cite{MR0159424}.
In this case, the coefficients
$(b_j,c_j)_{j=1}^s$
correspond to the
quadrature formula
of order $p$ such that
$\int_{0}^{1} \pi(\tau)\, \der \tau
=
\sum_{j=1}^s b_j \pi(c_j)$
for polynomials $\pi(\tau)$ of degree $\leq p-1$.  
Moreover, the coefficients in $\vec c$ and $\vec A$, together with conditions for the stage order $q$,
imply that
at every stage $i$ the quadrature formula
$\int_{0}^{c_i} \pi(\tau)\, \der \tau
=
\sum_{j=1}^s a_{ij} \pi (c_j)$
holds for polynomials $\pi(\tau)$ of degree $\leq q-1$.
Depending on the quadrature formula considered,
such as
Gauss, Radau or Lobatto, different families of
implicit Runge--Kutta methods can be constructed
(for more details, see \cite{Hairer96} \S~IV.5).

In this work we consider the family of RadauIIA
methods introduced by Ehle \cite{Ehle69},
based on \cite{Butcher64},
that consider
Radau quadrature formulas \cite{RadauC}
such that $p=2s-1$ and $q=s$.
These are
$A$-- and $L$--stable schemes that are
stiffly accurate methods according to
(\ref{eq:stiffly_acc}).
In particular we consider the third and fifth order RadauIIA schemes
referred as Radau3 and Radau5
(Butcher tableaux in Appendix~\ref{app:Butcher_tab}).
Note that, even though Gauss methods attain a maximum order
of $p=2s$ \cite{MR0159424,Ehle68},
they are neither
stiffly accurate nor $L$--stable schemes,
which are both important properties for stiff problems.
Approximations of lower order are obtained
with Lobatto methods satisfying
$p=2s-2$ \cite{MR0159424,Ehle68,Chipman71,Axelsson72}.
In particular the collocation
methods with $p=2s-2$ and $q=s$,
known as the LobattoIIIA methods, yield
stiffly accurate schemes, but these are only $A$--stable.

\section{Numerical implementation}\label{sec:imple}
We now discuss some particular aspects
concerning the numerical implementation.
We consider the multiresolution procedure
presented in \cite{Duarte11_SISC}.
For the sake of completeness some 
key details of this particular implementation
will be first recalled,
while more details and references
can be found in \cite{Duarte_Phd}.

\subsection{Construction of multiresolution grids}
The adapted grid is composed of a set
of nested dyadic grids: $S_j$, $j=0,1,\ldots,J$,
from the coarsest to the finest.
They are generated by recursively refining 
a given cell depending on the
local regularity of the 
time--dependent variables,
measured by the details
at a given time.
These grids
are implemented in a 
multi--dimensional Cartesian 
finite volume framework.
Following \cite{Cohen03}
a centered polynomial interpolation of accuracy order
$\beta = 2r+1$ is implemented
for the prediction operator, 
computed with the $r$ nearest neighboring cells
in each direction;
the procedure is exact for polynomials of degree $2r$.
Here we will only consider the case
$\beta=3$ with one neighboring cell per direction ($r=1$)
including the diagonals in multidimensional configurations.
For instance, in the one--dimensional case (\ref{eq3:dyadic_1D})
the latter is given by 
\begin{equation*}\label{eq3:polynomial_dyadic1D_order3}
\widehat{u}_{j+1,2k} = u_{j,k} + \frac{1}{8}(u_{j,k-1} - u_{j,k+1}), \qquad
\widehat{u}_{j+1,2k+1} = uf_{j,k} + \frac{1}{8}(u_{j,k+1} - u_{j,k-1}).
\end{equation*}
Higher order formulae can be found in \cite{muller2003}, 
while extension to multi-dimensional 
Cartesian grids is easily obtained by a
tensorial product of the one-dimensional operator 
\cite{bihari1997,Roussel03}.
In general the interpolation stencil is given
by $(2r+1)^d$ cells.

Data compression is achieved
by means of thresholding 
by discarding the cells whose details
are below a given tolerance
(see (\ref{eq3:MR_Lambda}))
and thus defining a compressed set $\Lambda$
with the remaining cells.
However, not all cells can be eliminated as this would prevent one 
from performing the multiresolution inter--grid operations.
In particular all cells in the prediction interpolation stencils 
must be always available.
Consequently, cells are gathered in a {\it graded tree} $\Lambda_\epsilon$,
instead of $\Lambda$,
that is, a data structure that satisfies the aforementioned conditions
(see \cite{Cohen03} for more details). 
Notice that $\Lambda \subset \Lambda_\epsilon$ and
error estimates like (\ref{eq3:adap_error_eps})
follow straightforwardly.
Nevertheless, for the ease of reading we will keep the 
notation $\Lambda$ in the following to 
refer to a graded tree.
A graded tree--structure is hence used
to represent data in the computer memory
(see also \cite{Roussel03}).
Recalling the standard tree--structure terminology:
if $\Omega_{\mu} \subset \Omega_{\lambda}$ 
belonging to consecutive grids, 
we say that $\Omega_\mu$ is a \textit{child} of 
$\Omega_\lambda$ and that $\Omega_\lambda$ is the \textit{parent} of $\Omega_\mu$. 
We thus define the \textit{leaves} 
$\leaf(\Lambda)$ of a 
\textit{tree} $\Lambda$ as the set of cells $\Omega_{\lambda}$,
$\lambda \in \leaf(\Lambda)$,
such that 
$\Omega_{\lambda}$ has no children in $\Lambda$.
Depending on the size of the computational domain, more 
graded trees may be needed.
Therefore,  cells are distributed in 
$N_{\rm R}$ graded trees $\Lambda_r$, $r=1,\ldots,N_{\rm R}$,
where $N_{\rm R}:=N_{{\rm R}x}N_{{\rm R}y}N_{{\rm R}z}$, 
and
$N_{{\rm R}x}$, $N_{{\rm R}y}$, and
$N_{{\rm R}z}$ stand for the number of graded trees
or {\it roots} per direction.
The adapted grid is thus given by sets
$\leaf(\Lambda_r)$, $r=1,\ldots,N_{\rm R}$,
with a total number of cells:
$N_\leaf= \sum_{r=1}^{N_{\rm R}}\#(\leaf(\Lambda_r))$.
If no adaptation is required, then the maximum number of cells will be
$N_{\rm L}=\#(S_J)=N_{{\rm R}x}N_{{\rm R}y}N_{{\rm R}z}2^{dJ}$\ie
the size of the finest grid.

Input parameters for the multiresolution implementation are: 
the maximum grid--level
$J$ corresponding to the finest spatial discretization;
the number of roots per direction
$N_{{\rm R}x}$, $N_{{\rm R}y}$, and
$N_{{\rm R}z}$;
and the threshold parameter
$\tolMR$, which defines the numerical
accuracy of the compressed representations following
(\ref{eq3:adap_error_eps}).

\subsection{Numerical function evaluations}

Introducing the set
${\rm I}^n_\leaf := \{1, 2, \ldots, N^n_\leaf\}$,
where $N^n_\leaf$ stands for the number of leaves
during time $t \in [t_n,t_n+\step_n]$,
we define a bijective function $h_n:D(h_n)\to{\rm I}^n_\leaf$,
with
\begin{equation*}
D(h_n):= \bigcup_{r=1}^{N_{\rm R}}  \leaf(\Lambda^n_r).
\end{equation*}
The cells 
$(\Omega_{\lambda})_{h_n(\lambda)\in {\rm I}^n_\leaf }$
correspond then to the adapted grid
during the current timestep $\step_n$, defined by 
the leaves of the tree representation.
Considering again $m=1$,
the solution of the semi--discrete problem (\ref{eq:gen_disc_prob})
for $t\in [t_n,t_n+\step_n]$
is similarly defined as $\disc U (t) = (u_{\lambda})_{h_n(\lambda)\in {\rm I}^n_\leaf }$,
where $u_{\lambda}$ stands for the
cell--average of variable $u(t,\x)$
in $\Omega _{\lambda}$ according to (\ref{eq3:average_finite_vol}).
The discrete function $\disc F (\disc U)$ 
in (\ref{eq:gen_disc_prob})
can be thus defined
as $\disc F (\disc U)= (F_{\lambda}(\disc U))_{h_n(\lambda)\in {\rm I}^n_\leaf }$,
where $F_{\lambda}(\disc U)$ can be further decomposed into
$\varPhi_{\lambda}(\disc U)$ and $\omega_{\lambda}(\disc U)$,
coming from the discretization of differential operators
and source terms, respectively.
In particular for a second order spatial discretization, 
considered in this work, the local source term $\omega_{\lambda}(\disc U)$
becomes $\omega(u_{\lambda})$, that is, it is computed using the local 
cell--average values.
During timestep $\step_n$, 
the time--dependent problem (\ref{eq:gen_disc_prob})
can be thus written at each cell $\Omega _{\lambda}$
of the adapted grid as
\begin{equation}\label{eq:disc_prob_local}
\der_t u_{\lambda} = 
F_{\lambda}(\disc U)= 
\varPhi_{\lambda}(\disc U) + \omega(u_{\lambda}), 
\qquad
t \in [t_n,t_n+\step_n],\, h_n(\lambda)\in {\rm I}^n_\leaf,
\end{equation}
\begin{equation}\label{eq:def_flux}
\varPhi_{\lambda}(\disc U) :=
|\Omega_{\lambda}|^{-1} \sum_{\mu}
|\Gamma_{\lambda,\mu}| \varPhi_{\lambda,\mu}, 
\end{equation}
where the latter sum is made over all $\mu \neq \lambda$
such that the interface 
$\Gamma_{\lambda,\mu} := \overline{\Omega_{\lambda}} \cap 
\overline{\Omega_{\mu}}$ is not trivial\ie
over all the neighboring cells of $\Omega_{\lambda}$;
and $\varPhi_{\lambda,\mu}$ accounts for the flux across each interface.
In the simplest 
(low order in space) schemes,
the flux $\varPhi_{\lambda,\mu}$ is typically a function of 
$u_{\lambda}$ and $u_{\mu}$ only, 
while higher order schemes require considering additional
cells.

Without loss of generality,
let us denote by $R_{\varPhi}[\lambda]$ 
the stencil required to compute fluxes 
associated with cell $\Omega_{\lambda}$.
Here, we consider flux computation 
schemes such that
all cells in $R_{\varPhi}[\lambda]$ belong to the same
grid, that is, 
fluxes are computed on a locally uniform mesh.
Problem (\ref{eq:disc_prob_local}) can be thus rewritten as
\begin{equation}\label{eq:disc_prob_local_flux}
\der_t u_{\lambda} = F_{\lambda}\left(\left(u_{\lambda}\right)_{\lambda \in R_{\varPhi}[\lambda]}\right), 
\qquad
t \in [t_n,t_n+\step_n],\, h_n(\lambda)\in {\rm I}^n_\leaf.
\end{equation}
The numerical integration of problem (\ref{eq:gen_disc_prob})
then involves evaluating function 
$F_{\lambda}$ in (\ref{eq:disc_prob_local_flux})
for the $N^n_\leaf$ current cells.
Moreover,
for a given interface $\Gamma_{\lambda,\mu}$
the following conservation property holds
in a finite volume flux representation:
$\varPhi_{\lambda,\mu} +  \varPhi_{\mu,\lambda} = 0$.
Computing $\varPhi_{\lambda,\mu}$ for $\Omega_\lambda$
amounts to evaluating also $\varPhi_{\mu,\lambda}$ for the neighboring
cell $\Omega_\mu$. 
Let us denote $\varPhi_{\lambda,\mu}^+$ as the right flux for 
$\Omega_\lambda$ and $\varPhi_{\mu,\lambda}^-$ as the left flux
for $\Omega_\mu$, along the direction normal to $\Gamma_{\lambda,\mu}$.
Similarly, $R_{\varPhi}^+[\lambda]$ stands for the stencil required to
compute $\varPhi_{\lambda,\mu}^+$ and, naturally, 
$R^-_{\varPhi}[\mu] \equiv R^+_{\varPhi}[\lambda]$;
we thus have that 
$\varPhi_{\mu,\lambda}^-= -\varPhi_{\lambda,\mu}^+$.
This property is thus exploited 
to save computations as
fluxes are computed only once at each interface.
The locally uniform grids are then defined by the stencil $R_{\varPhi}^+[\lambda]$
enclosing the current leaf $\Omega_{\lambda}$.
Ghost cells, computed according to the inter--grid prediction operation,
are used whenever one cell in the current stencil is missing.
These ghost cells are also added 
to the adapted grid at interfaces between cells of different sizes
in order to compute numerical fluxes at the
highest grid--level between two neighboring cells \cite{Roussel03}.

Notice that function $h_n$ is in practice used for indexation of leaves,
identifying them regardless of their geometric layout.
This is particularly useful to organize the computation of the entries
of the Jacobian, as shown in Appendix~\ref{app:Jac},
and the linear system.
All matrices are stored using a standard CSR (Compressed
Sparse Row) format for sparse matrices.

\subsection{Newton method and linear solver}
The simplified Newton method to 
solve the nonlinear system (\ref{eq2:nonlinear_sys})
considers the linear system (\ref{eq2:Newton_simp})
for fully IRK schemes like the RadauIIA methods.
System (\ref{eq2:Newton_simp}) is recast as
\begin{equation}\label{eq2:Newton_simp_dt}
(\step^{-1}\Ibf_{s\times m \times N} - \disc A \disc J) \delta \disc Z^k 
= -\step^{-1}\disc Z^k+ \disc A \calbf F(\disc Z^k),
\end{equation}
mainly to avoid updating 
$\step \disc A \disc J$ when time step changes are required.
Defining an accuracy tolerance 
$\tolnew$, we consider the following stopping criterion
for the iterative process:
\begin{equation}\label{eq:stop_newton}
\|\delta \disc Z^k \|_2 \leq \tolnew.
\end{equation}
Additionally, we define a convergence rate for the Newton solver as
\begin{equation}\label{eq:conv_rate}
\Theta_k = \frac{\|\delta \disc Z^k \|_2}{\|\delta \disc Z^{k-1} \|_2},
 \quad k\geq 1.
\end{equation}
For the first iteration we set 
$\Theta_0 = \|\delta \disc Z^0 \|_2/(2\max \disc U_0)$.
We also define a maximum number of Newton 
iterations $\knew$, and inspired by \cite{Hairer96},
computations are interrupted and restarted with a 
halved timestep in (\ref{eq2:Newton_simp_dt}),
if any of the following happens:
\begin{itemize}
 \item there is a $k$ such that $\Theta_k \geq 1$;
 \item for some $k$, we have that
 \begin{equation}\label{eq:err_new_max}
  (\Theta_k)^{\knew-k-1}\|\delta \disc Z^k \|_2 \geq \tolnew,
 \end{equation}
 where the left--hand side in (\ref{eq:err_new_max})
 is a rough estimate of $\|\delta \disc Z^{\knew -1} \|_2$;
 \item $\knew$ iterations have been performed and 
 $\|\delta \disc Z^{\knew-1} \|_2 > \tolnew$.
\end{itemize}
Notice that if the timestep is halved, only diagonal entries
in $[\step^{-1}\Ibf_{s\times m \times N} - \disc A \disc J]$
need to be modified;
however, the resulting new matrix must be again factorized.

In this work we have implemented the iterative 
GMRES method \cite{GMRES} to solve the linear system
(\ref{eq2:Newton_simp_dt}),
with right-preconditioning based on an ILUT factorization \cite{ILUT}.
Considering a fixed Jacobian has the advantage that the 
factorization and preconditioning
of matrix $[\step^{-1}\Ibf_{s\times m \times N} - \disc A \disc J]$
needs to be performed only once, 
unless computations are restarted with a halved timestep.
Notice that this is a purely algebraic problem that has completely 
lost any reminiscence of its original geometric layout,
meaning that it is independent of the adapted grid generation 
or any other grid--related data structure or geometric consideration.
Consequently, any linear solver could be used as a {\it black box}
solver provided that it only needs the matrix entries and the
right--hand side array as inputs.
For an iterative linear solver like GMRES 
we define another accuracy tolerance,
$\tollin$, as stopping criterion.
This tolerance is chosen such that
$\tollin=\kappa \tolnew$, with $\kappa \leq 1$.
In this work we consider, for instance, $\kappa = 10^{-2}$,
unless noted otherwise.
If the linear solver is taking too many iterations to converge,
noted henceforth as $\klinjac$ iterations, we
update the Jacobians in the Newton method.
Re--factorization and preconditioning would be necessary in this case.

The same ideas apply for the numerical implementation 
of the SDIRK schemes,
considering at each stage system
\begin{align}\label{eq2:Newton_simp_SDIRK_dt}
((\step\gamma)^{-1}\Ibf_{m \times N} - \disc J_0) \delta \disc z_i^k 
= & -(\step\gamma)^{-1}\disc z_i^k+ \disc F(t_0+c_i\step,\disc U_0+ \disc z^k_i ) \nonumber \\
& + \sum_{j=1}^{i-1} \frac{a_{ij}}{\gamma}\disc F(t_0+c_j\step,\disc U_0+ \disc z_j),
\end{align}
instead of (\ref{eq2:Newton_simp_SDIRK}).
In principle
the same block--matrix needs to be factorized
for all Newton iterations at all stages.
The stopping criterion (\ref{eq:stop_newton}),
the convergence rate (\ref{eq:conv_rate}), as well as
the conditions for halving the time step are all applied
stage--wise, that is, with 
$\|\delta \disc z^{k}_i \|_2$ instead of $\|\delta \disc Z^{k} \|_2$.
Similarly, the Jacobian is updated 
in (\ref{eq2:Newton_simp_SDIRK_dt})
after $\klinjac$ iterations.

\subsection{Time--stepping strategy}
Since we consider only $A$--stable IRK schemes in this work,
we are uninhibited by stability issues in the choice of the timestep,
which can be based solely on accuracy requirements.
For some kinds of problems, 
a constant time step might be sufficient 
to capture the problem dynamics.
However, more generally, 
an adaptive time--stepping could
be considered 
in order to enhance the computational efficiency.
In either case, the main goal is to define
a time step $\step$ such that the local
error satisfies
\begin{equation}\label{eq2:tolerance_err}
\|\disc U(t_0+\step) - \disc U_1\|_2 = C \step^{p+1} \leq \toltime,
\end{equation}
where $\toltime$ is the desired accuracy tolerance for the
$p$--th order IRK scheme.
The advantage of higher order methods is that they can 
satisfy (\ref{eq2:tolerance_err}) with larger time steps
than those achievable with conventional low order methods.
A standard approach to time step control is based on 
numerically approximating
the exact local error
in (\ref{eq2:tolerance_err}), by 
considering a solution $\hat{\disc U}_1$ computed by a 
lower order method of order $\hat{p}<p$ (see, for instance, \cite{Hairer87}).
In this way we use the computations at the $n$--th step to predict
the local error at the next step,
\begin{equation}\label{eq:err_estimate}
 err = \| \disc U_n - \hat{\disc U}_n \|_2 \approx  \tilde{C}_n \step_n^{\hat{p}+1},
\end{equation}
which defines a new time step,
\begin{equation}\label{eq2:time_stepping1_n}
\step_{\new} = \step_n \ds \left(\frac{\toltime}{err}\right)^{1/\hat{p}+1},
\end{equation}
by assuming that 
$ \toltime \approx \tilde{C}_{n+1} \step_{\new}^{\hat{p}+1}$
with $\tilde{C}_{n+1} \approx \tilde{C}_n $.
The next time step $\step_{n+1}$ will be based on $\step_{\new}$
if the current approximation error satisfies
$err \leq \toltime$.
Otherwise,
the current $n$--th solution will be
disregarded, and 
the same $n$--th step
will be integrated again with 
$\step_{\new}$ instead of $\step_n$.
The lower order approximations are defined in Appendix~\ref{app:embeddedRK}.

Inspired by \cite{Hairer96}, we define a safety factor $\nu_k$
that depends on the 
current Newton iteration
$k$, the current linear solver iteration $\klin$,
and the
maximum number of Newton iterations
$\knew$,
as follows
\begin{equation}\label{eq:def_nu}
 \nu_k = \nu \times \frac{2\knew +1}{2\knew + \max(k,0.5\klin)}, 
\end{equation}
where $\nu>0$ is a standard safety factor close to 1.
Here we typically consider $\nu = 0.9$.
For the SDIRK schemes, where more than one Newton solve
is required per time step, $k$ and $\klin$ in (\ref{eq:def_nu})
stand, respectively, for the maximum number of Newton and linear 
iterations performed
within a given time step.
The time step $\step_{n+1}$ is thus defined as
\begin{equation}\label{eq:timestepping_dyn}
 \step_{n+1} = \min(\nu_k \step_{\new}, \alpha \step_{n}),
\end{equation}
where $\alpha>1$ limits the variation of successive time steps.
Here we consider in general $\alpha = 1.5$.
However, if the computations were to be performed with a constant 
time step $\step$, we would consider the following time--stepping procedure
\begin{equation}\label{eq:timestepping_const}
 \step_{n+1} = \min(\alpha \nu_k \step_{n}, \step),
\end{equation}
which allows modifications on the chosen time step based on the
performance of the Newton and linear solvers.
In general the initial time step $\step_0$ should be set 
sufficiently small to account for potentially fast transients.

The numerical accuracy of the time integration is defined
by the user--provided tolerance parameter, $\toltime$.
The tolerance parameter for the Newton solver 
is set to a lower value: $\tolnew = \kappa \toltime$,
with $\kappa < 1$.
In this way errors coming from both the Newton 
and linear solvers should remain
smaller than those caused by the IRK scheme.

\section{Numerical illustrations}\label{sec:num_res}

We investigate the computational performance of the numerical
strategy for three problems modeled by time--dependent
stiff PDEs.
In this work all the simulations were run on a standard laptop
with an Intel Core i3 @ $2.27$ GHz processor and a memory capacity of $1.8$ GB.

\subsection{The Belousov--Zhabotinski reaction}\label{subsec:BZ}

Let us consider the numerical approximation of a model for
the Belousov--Zha\-bo\-tins\-ki (BZ) reaction, a catalyzed oxidation of an organic species by acid
bromated ion
(see \cite{Epstein98} for more details and illustrations).
The present mathematical formulation \cite{Field72,Scott94}
takes into account three species: hypobromous acid $\mathrm{HBrO_2}$,
bromide ions $\mathrm{Br^-}$, and cerium (IV).
Denoting by $a=[\mathrm{Ce(IV)}]$, $b=[\mathrm{HBrO_2}]$, and $c=[\mathrm{Br^-}]$,
we obtain a very stiff system of three PDEs given by
\begin{equation} \label{eq4:bz_eq_3var_diff}
\left.
\begin{array}{l}
\partial_t a - D_a\, \partial^2 _\x a = \ds \frac{1}{\mu}(-qa-ab+fc),\\[1.75ex]
\partial_t b - D_b\, \partial^2 _\x b = \ds \frac{1}{\varepsilon}\left(qa-ab+b(1-b)\right),\\[1.75ex]
\partial_t c - D_c\, \partial^2 _\x c = b-c,
\end{array}
\right\}
\end{equation}
where $\x \in \R^d$,
with real, positive
parameters: $f$, small $q$, and small $\varepsilon$ and
$\mu$, such that  $\mu \ll \varepsilon \ll 1$.
In this study: $\varepsilon = 10^{-2}$,
$\mu = 10^{-5}$, $f=1.6$, $q=2\times 10^{-3}$;
with diffusion coefficients:
$D_a=2.5\times 10^{-3}$, $D_b=2.5\times 10^{-3}$, and $D_c=1.5\times 10^{-3}$.
The dynamical system associated with this problem
models reactive, excitable media with a
large time scale spectrum (see \cite{Scott94}
for more details).
The spatial configuration with
the addition of diffusion involves propagating wavefronts
with steep spatial gradients;
in particular, two--dimensional spiral waves and three--dimensional
scroll waves \cite{Duarte11_SISC}.

\subsubsection{Numerical time integration errors}\label{subsubsec:BZ_1D}
\begin{figure}[!htb]
\begin{center}
\includegraphics[width=0.49\textwidth]{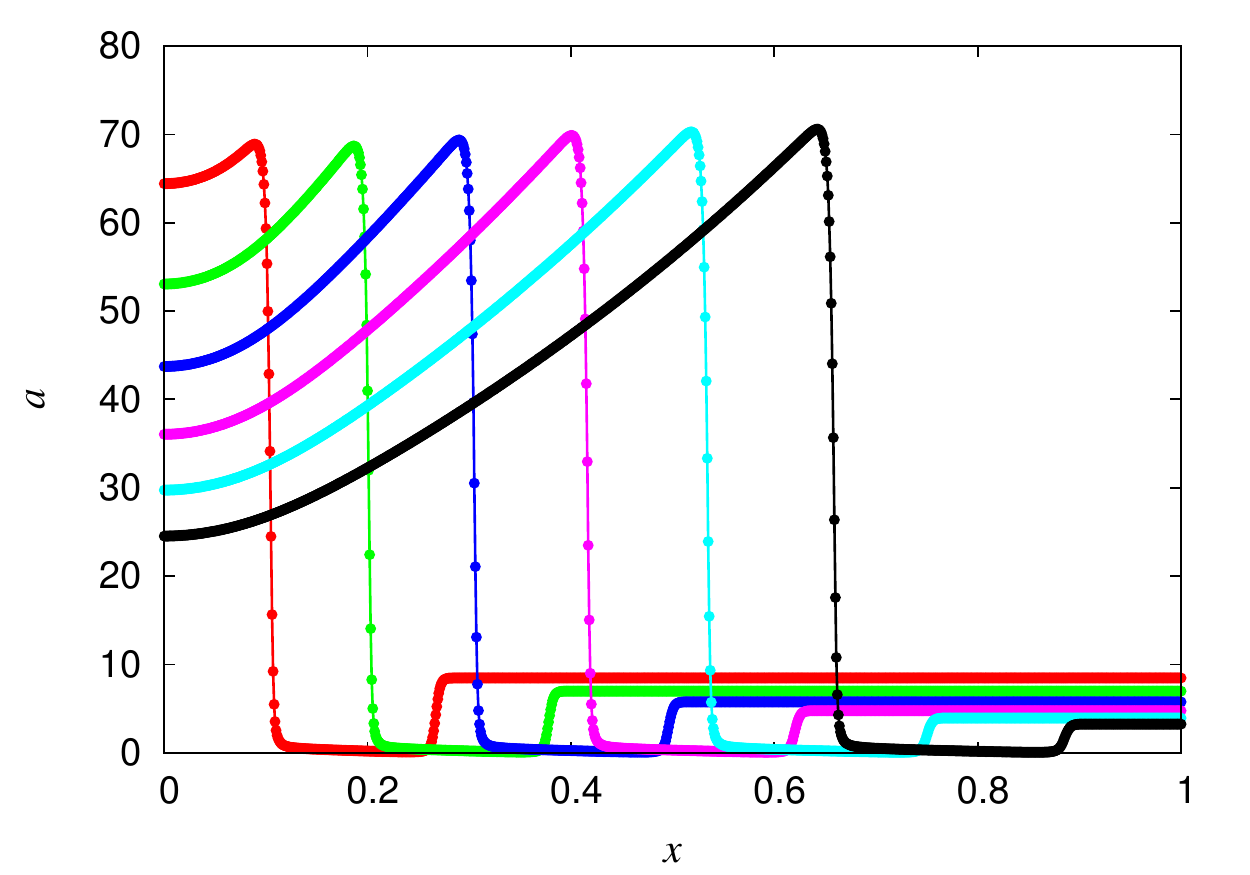}
\includegraphics[width=0.49\textwidth]{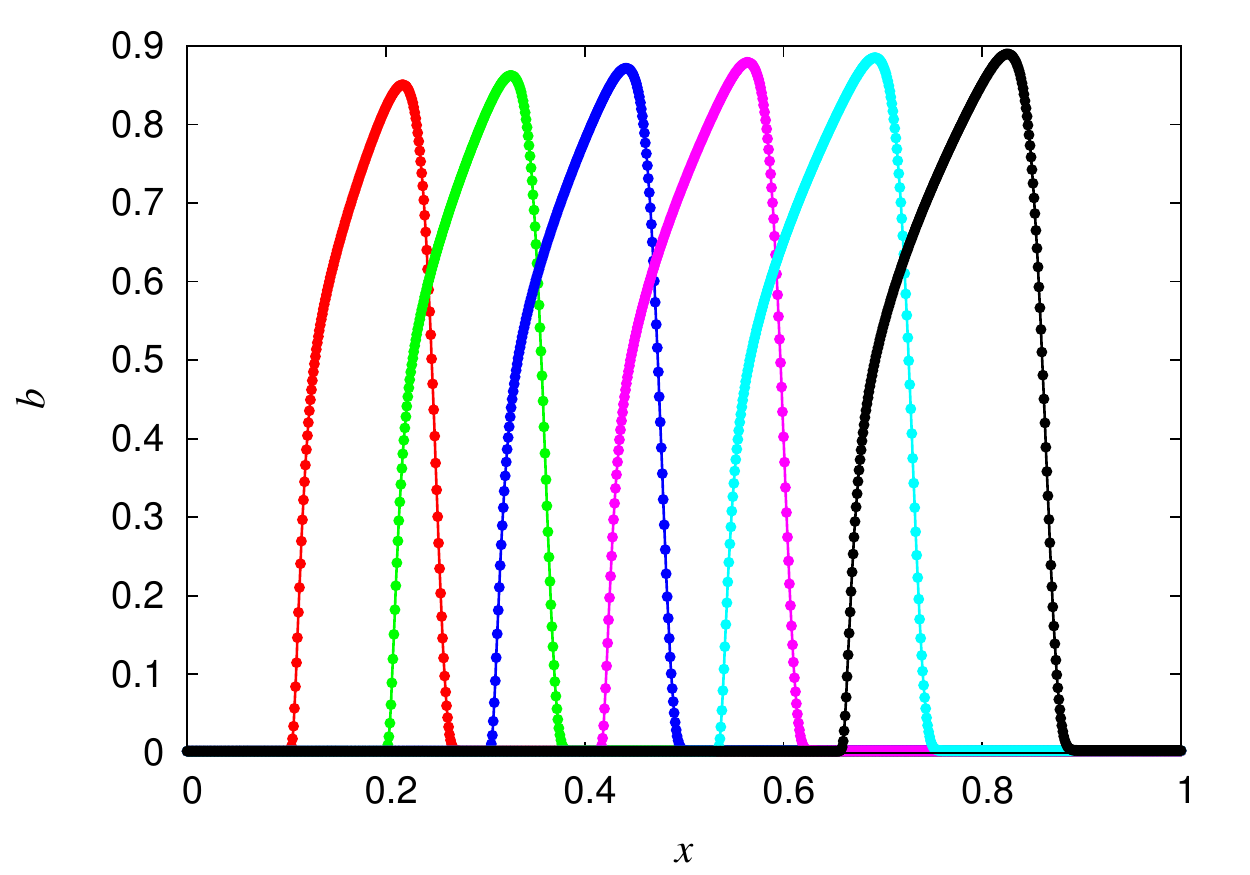}
\includegraphics[width=0.49\textwidth]{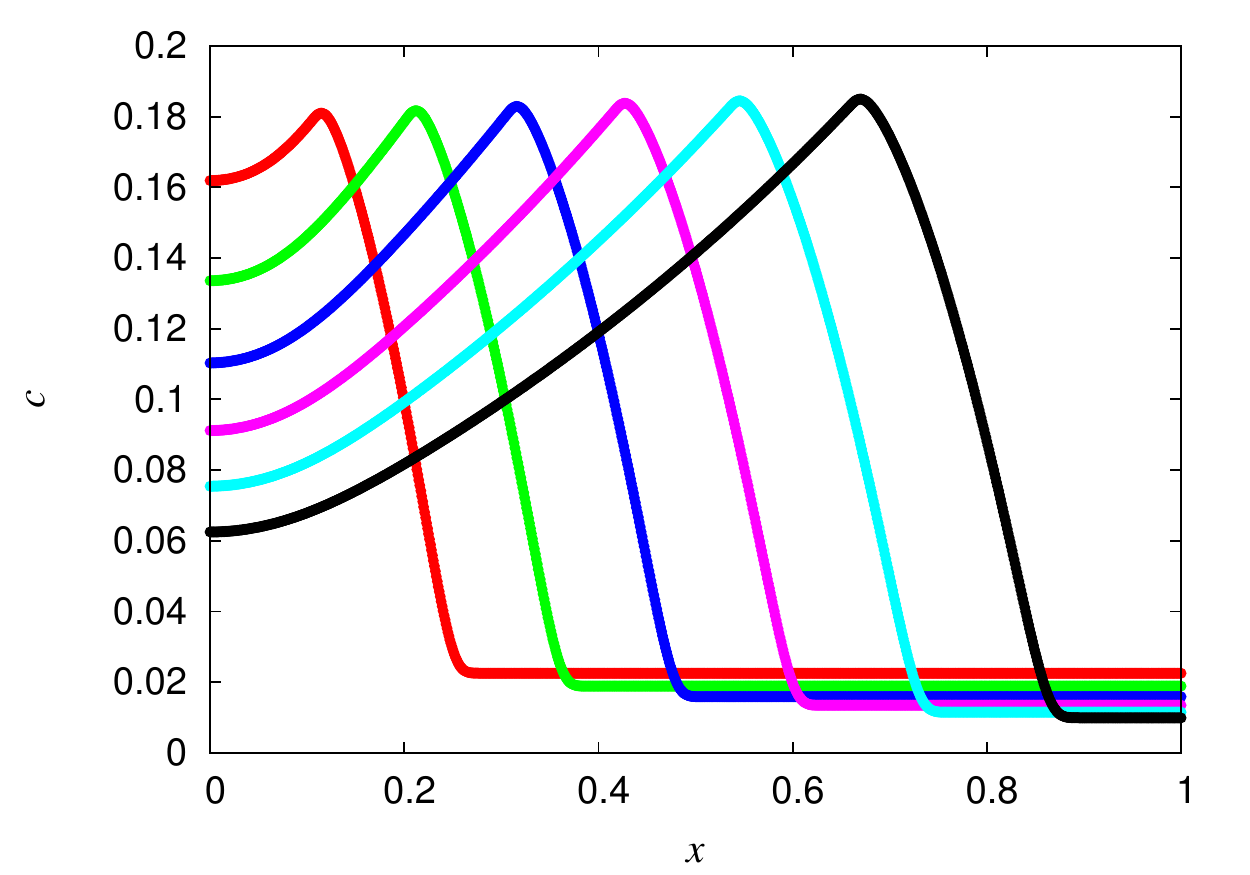}
\end{center}
\caption{One--dimensional BZ propagating waves for variables $a$ (top left),
$b$ (top right), and $c$ (bottom), at time intervals of $0.2$
within $[0,1]$ from left to right.}
\label{fig:sol_BZ_nx1001}
\end{figure}

We consider problem
(\ref{eq4:bz_eq_3var_diff}) in a one--dimensional configuration
with Neumann homogeneous boundary conditions,
discretized on a uniform grid of 1024 cells
over a space region of $[0,1]$.
A standard, second order, centered finite volumes scheme
is employed for the diffusion term.
No grid adaptation is considered here in order to assess 
only the numerical errors related to the time integration schemes.
To obtain an initial condition, we initialize the problem
with a discontinuous profile close to the left boundary;
we then integrate in time until the BZ wavefronts are fully developed.
Figure~\ref{fig:sol_BZ_nx1001} shows the time evolution of the
propagating waves for a time window of $[0,1]$.
In order to compute the local errors associated with the
implicit solvers here considered,
we define a reference solution for the
resulting semi--discrete problem. 
The latter is chosen here as
the solution obtained using the Radau5 scheme (\ref{eq:Radau5}),
computed with a fine tolerance: $\toltime = 10^{-14}$.
\begin{figure}[!htb]
\begin{center}
\includegraphics[width=0.49\textwidth]{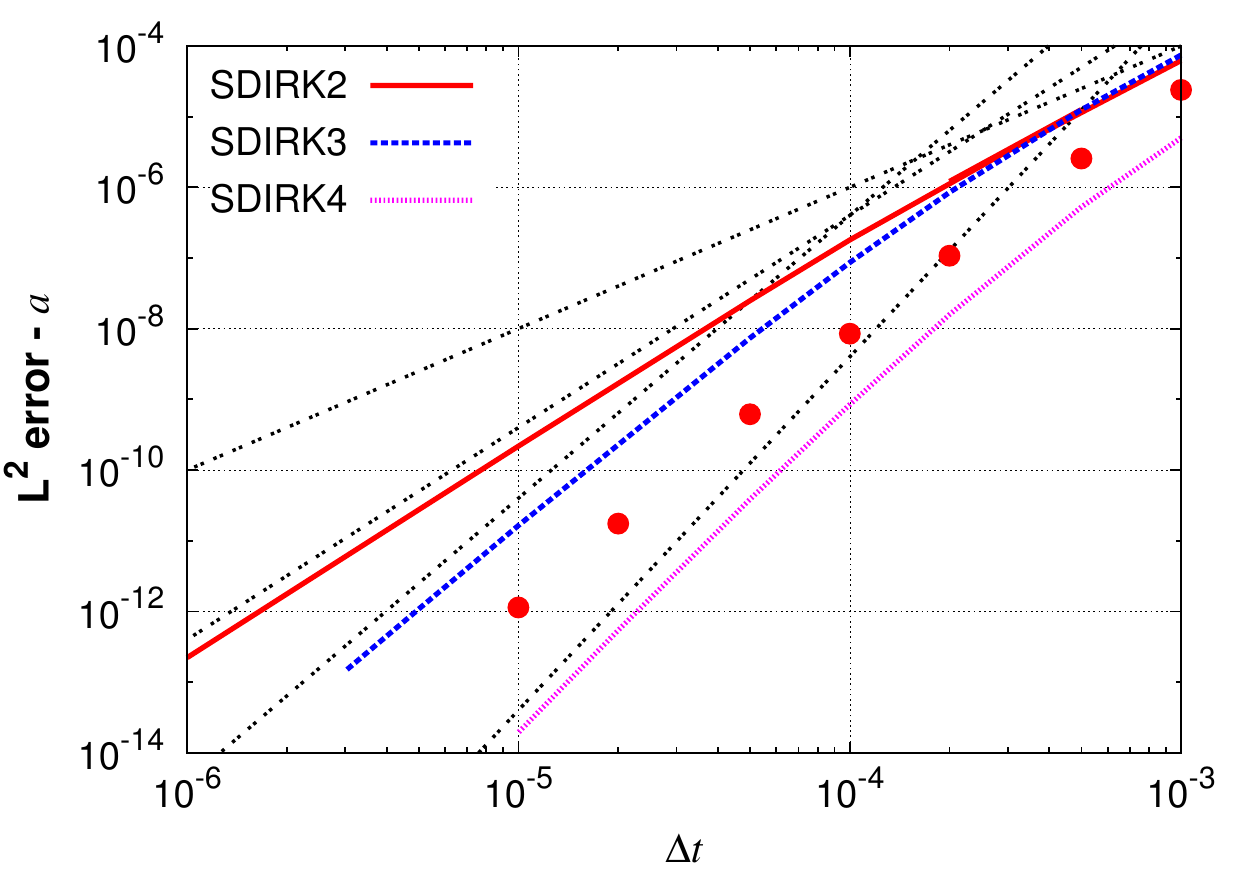}
\includegraphics[width=0.49\textwidth]{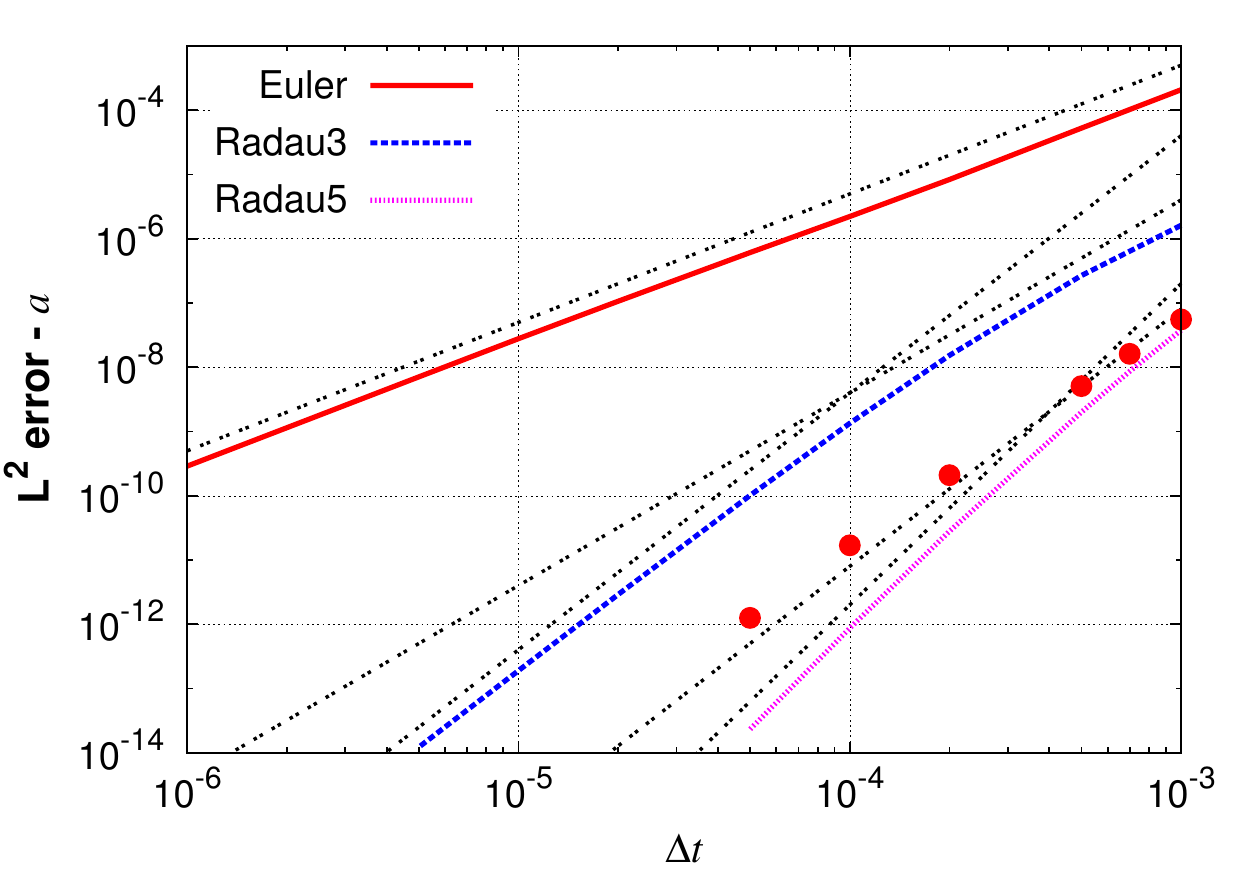}
\includegraphics[width=0.49\textwidth]{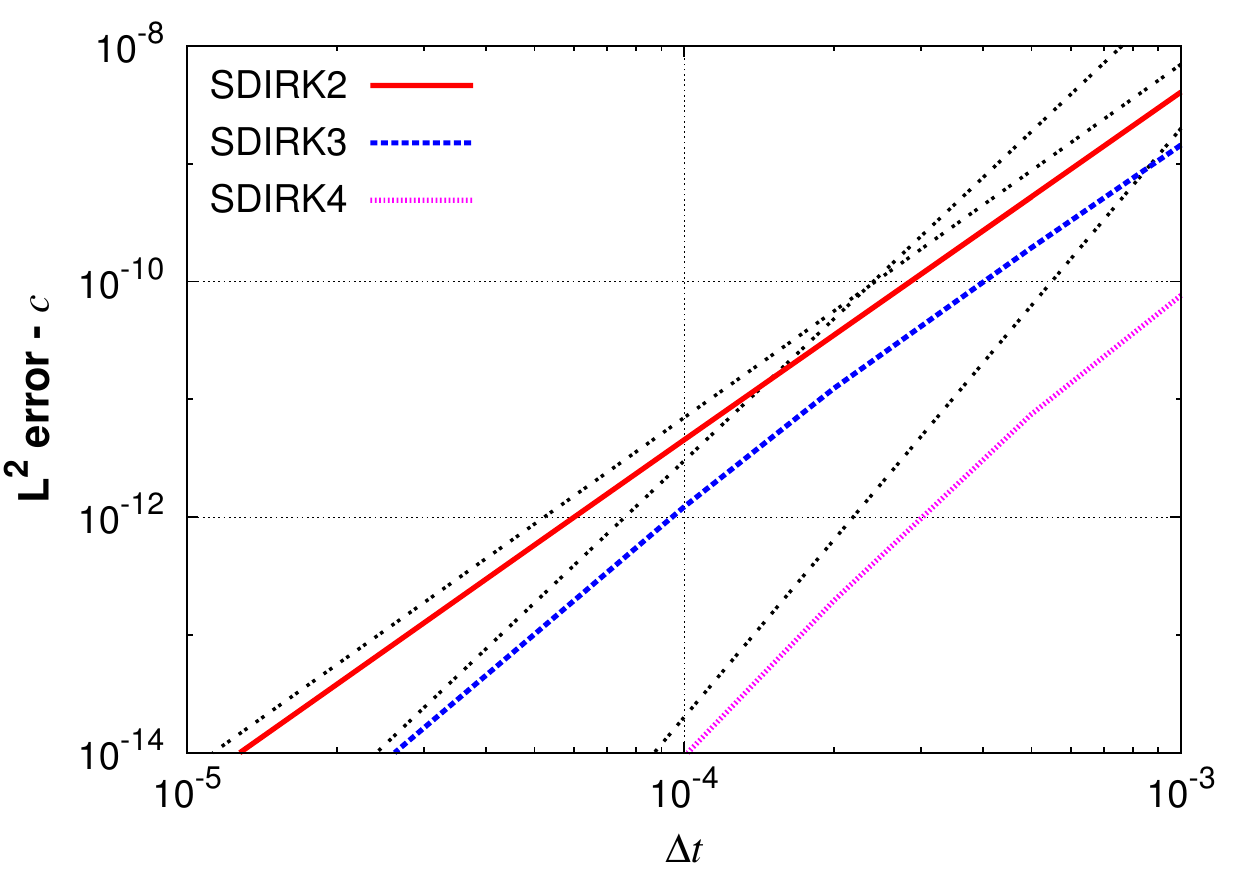}
\includegraphics[width=0.49\textwidth]{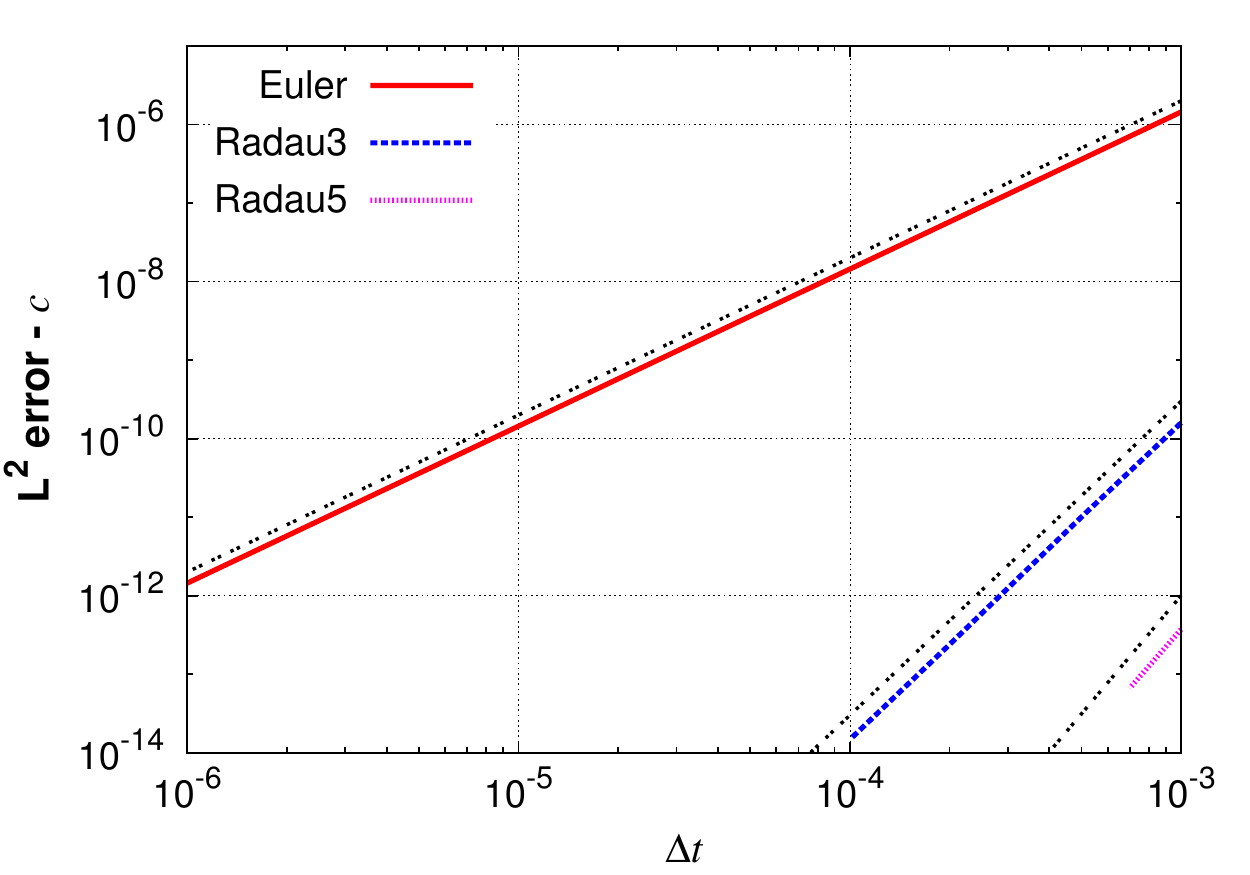}
\end{center}
\caption{Local $L^2$--errors for stiff 
and non--stiff components, respectively, 
$a$ (top) and  $c$ (bottom) using the SDIRK2,
SDIRK3 and SDIRK4 schemes (left), and the 
Euler, Radau3 and Radau5 ones (right).
Dashed lines of slopes 2 to 5 (top),
3 to 5 (bottom left), and
2, 4 and 5 (bottom right)
are also depicted.
Error estimates $err$ given by (\ref{eq:err_estimate})
are indicated with red bullets ($\color{red}{\bullet }$)
(top) for the SDIRK4 (left)
and Radau5 (right) schemes.}
\label{fig:local_order_BZ1D}
\end{figure}

Starting from the solution at $t=0.5$, Figure~\ref{fig:local_order_BZ1D} 
shows the local errors associated with each IRK scheme
for different time steps.
Both tolerances for the Newton and the linear solver are
set to $\tolnew = \tollin = 10^{-14}$ in these computations.
Notice that the stiffest variable, $a$, is directly subject to a time scale
given by the small parameter $\mu = 10^{-5}$.
We are thus in practice 
interested in time steps larger than $10^{-5}$.
Considering the stiff and non--stiff components of (\ref{eq4:bz_eq_3var_diff}), 
$a$ and $c$, respectively, we see the following
numerical behavior.
For the stiff variable (see Figure~\ref{fig:local_order_BZ1D} (top)),
local errors of $\Or(\step^{p+1})$ tend to $\Or(\step^{q+1})$ for 
relatively large time steps.
For the non--stiff variable (see Figure~\ref{fig:local_order_BZ1D} (bottom)),
the order reduction goes from $\Or(\step^{p+1})$ to $\Or(\step^{q+2})$.
Local errors for variable $b$ (not shown), 
stiffer than $c$, also behave as the ones
for the stiffest component, variable $a$.
These results are consistent with the classical, theoretical bounds derived
in \cite{Hairer88} for stiff ODEs in {\it singular perturbation} form,
that is, containing a small stiffness parameter given by $\mu$
in our case.
These results highlight the importance of the stage order for 
IRK schemes and stiff problems.
In this respect RadauIIA schemes perform better than SDIRK methods.
The same can be said with respect to stiffly accurate schemes
when comparing, for instance, SDIRK3 with Radau3;
more accurate results are obtained with the latter.
As a matter of fact, a well--known conclusion is that 
stiffly accurate schemes guarantee better accuracies
for stiff problems \cite{Prothero74,Alexander77,Hairer88}.
Figure~\ref{fig:local_order_BZ1D} also shows the 
error estimates $err$ given by (\ref{eq:err_estimate})
for both SDIRK4 and Radau5 schemes.
Notice that the actual local errors are bounded by $err$,
which in particular overestimates them since 
$err$ is computed using a third order, embedded scheme in 
both cases.
Finally, it has to be remarked that
higher order schemes perform better in terms of numerical
accuracy than low order ones like the first order Euler method,
even when order reduction appears and all methods show the 
same low order convergence.

\subsubsection{Performance comparison}
\begin{figure}[!htb]
\begin{center}
\includegraphics[width=0.495\textwidth]{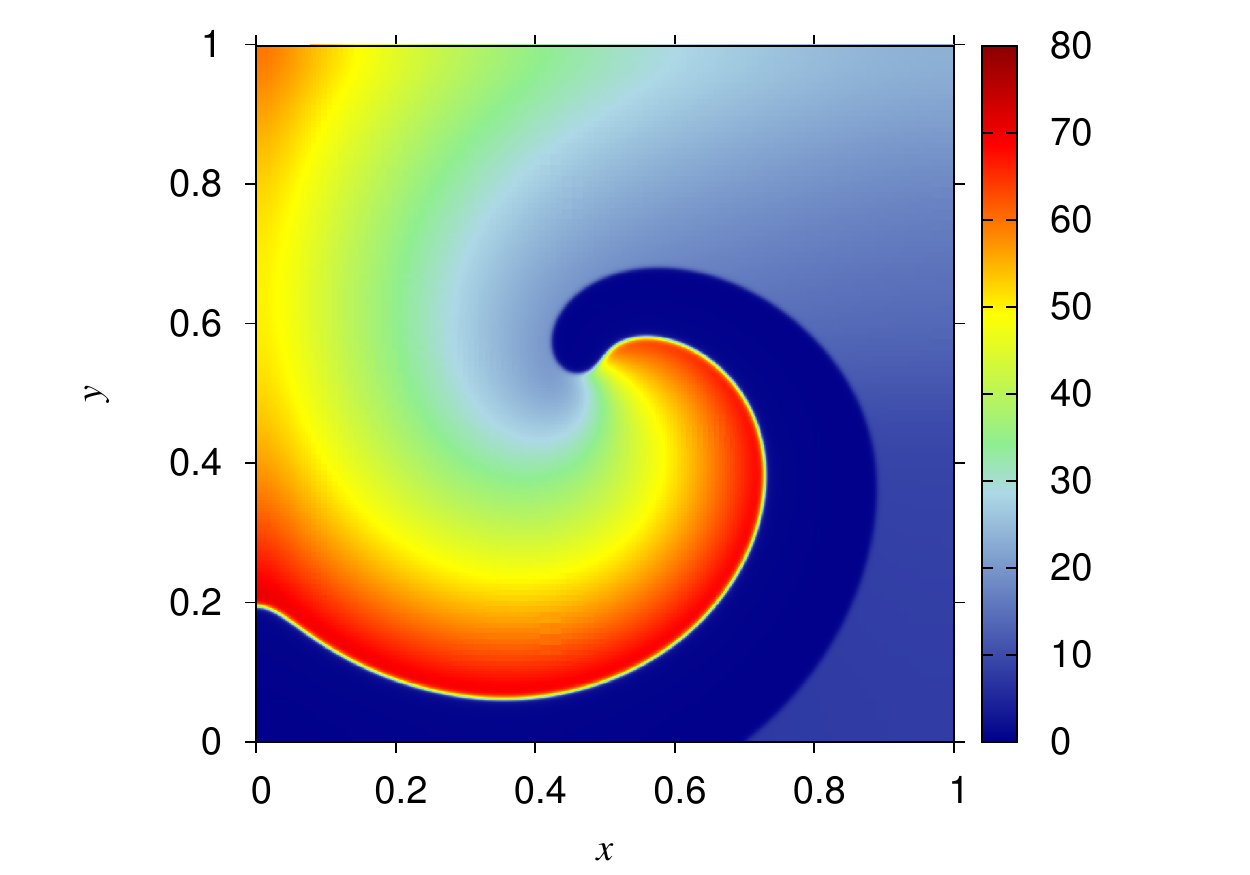}
\includegraphics[width=0.495\textwidth]{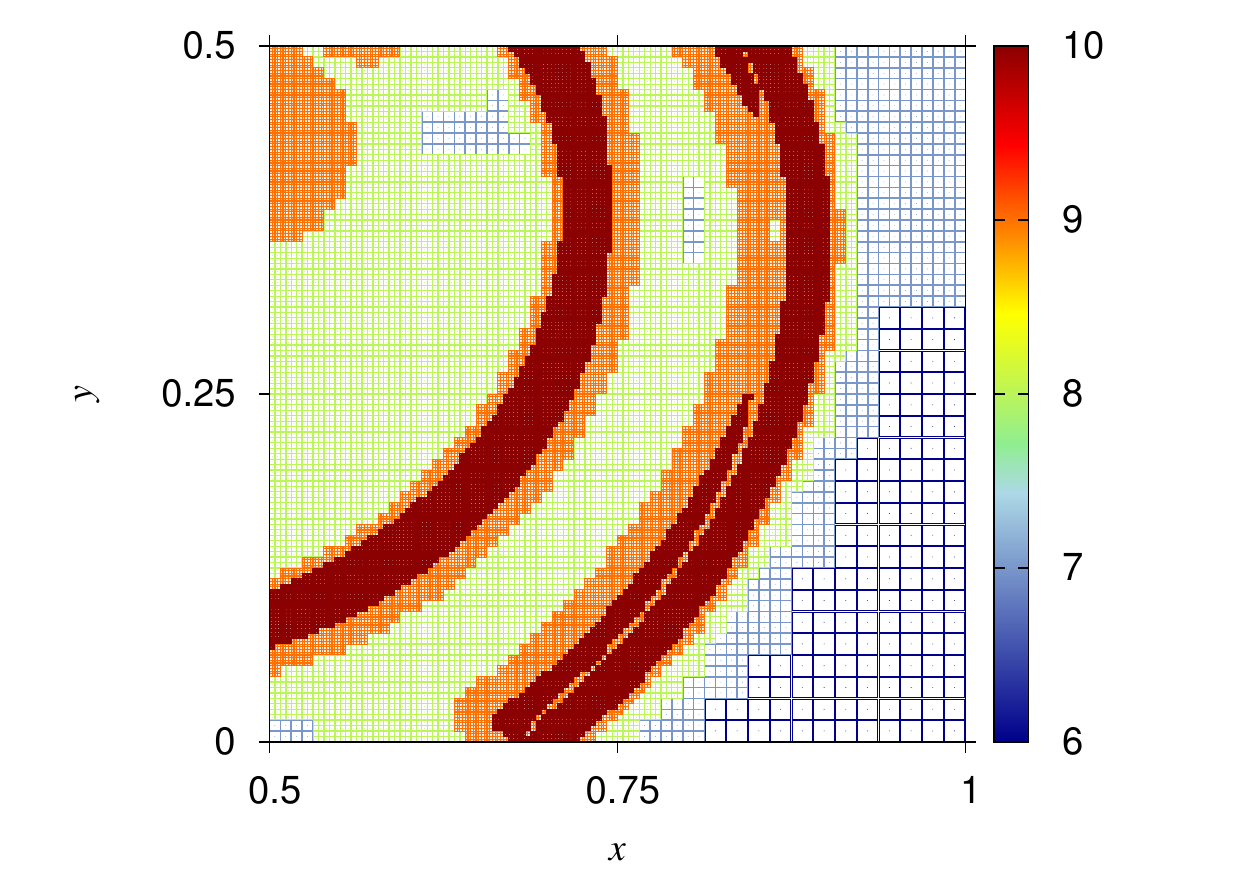}
\end{center}
\caption{Two--dimensional BZ propagating waves for variable $a$ at $t=2$ (left)
and the corresponding adapted grid for the zoomed region $[0.5,1]\times[0,0.5]$ (right).}
\label{fig:sol_BZ_2D}
\end{figure}

We now consider problem
(\ref{eq4:bz_eq_3var_diff}) in a two--dimensional configuration
with Neumann homogeneous boundary conditions,
using multiresolution analysis to adapt dynamically the 
spatial discretization grid.
For the multiresolution analysis the following input
parameters are considered:
number of roots per direction,
$N_{{\rm R}x}=N_{{\rm R}y}=1$; 
maximum grid--level, $J=10$; and
accuracy tolerance,
$\tolMR = 10^{-3}$.
The finest grid has a spatial resolution of 
$1024 \times 1024$ over a computational domain of $[0,1]\times [0,1]$.
We consider both SDIRK4 and Radau5 schemes with the following
parameters:
$\knew = 30$, $\klinjac = \knew$, and
$\kappa = 10^{-1}$, recalling that 
$\tollin = \kappa \tolnew = \kappa^2 \toltime$.
The initial solution is taken at $t=2$, when the spiral waves
are fully developed
(see Figure~\ref{fig:sol_BZ_2D}), and the PDEs are then
integrated until $t=2.01$.
(See \cite{Duarte11_SISC} for details on the initialization 
of this two--dimensional configuration.)
The data compression,
defined as the ratio in percentage
between the active and the finest grids,
is of about 15\%.
\begin{figure}[!htb]
\begin{center}
\includegraphics[width=0.49\textwidth]{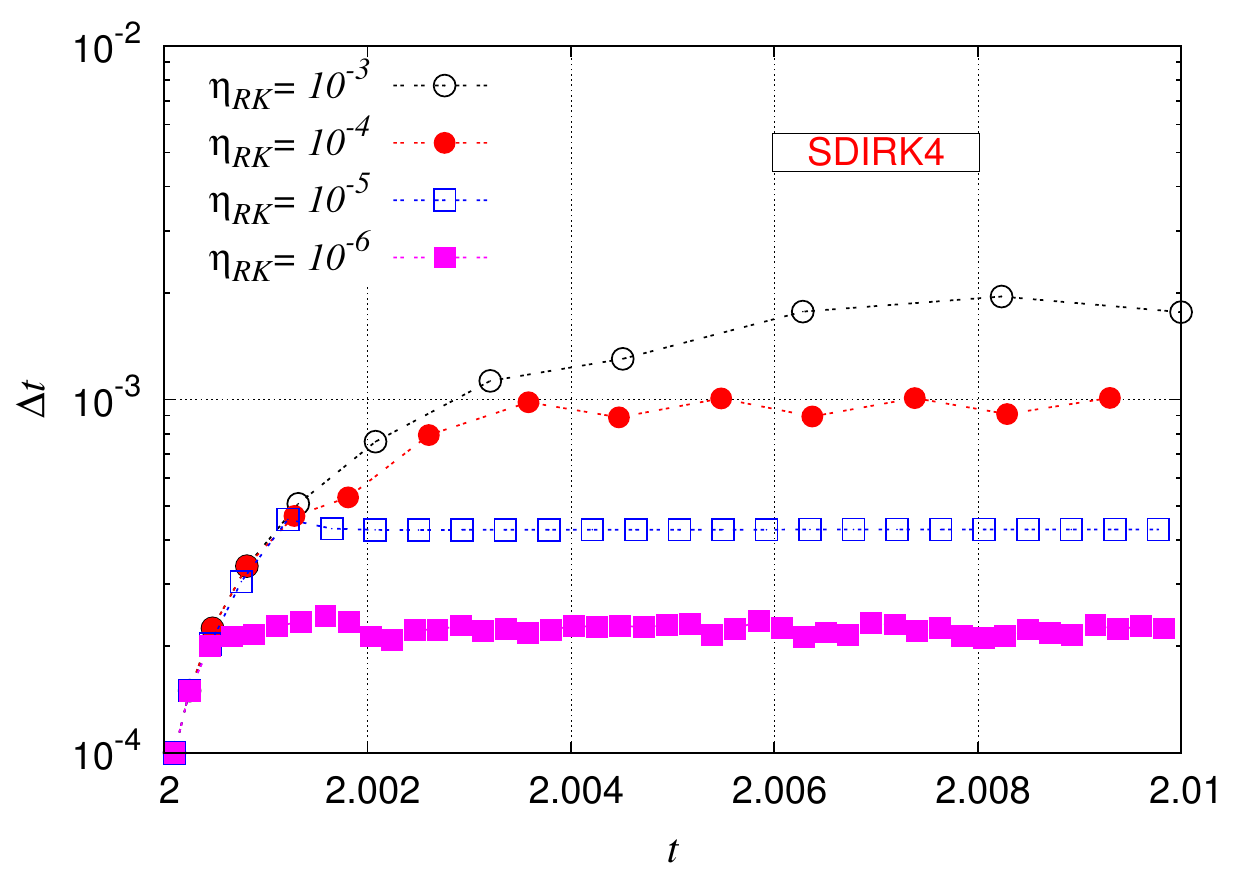}
\includegraphics[width=0.49\textwidth]{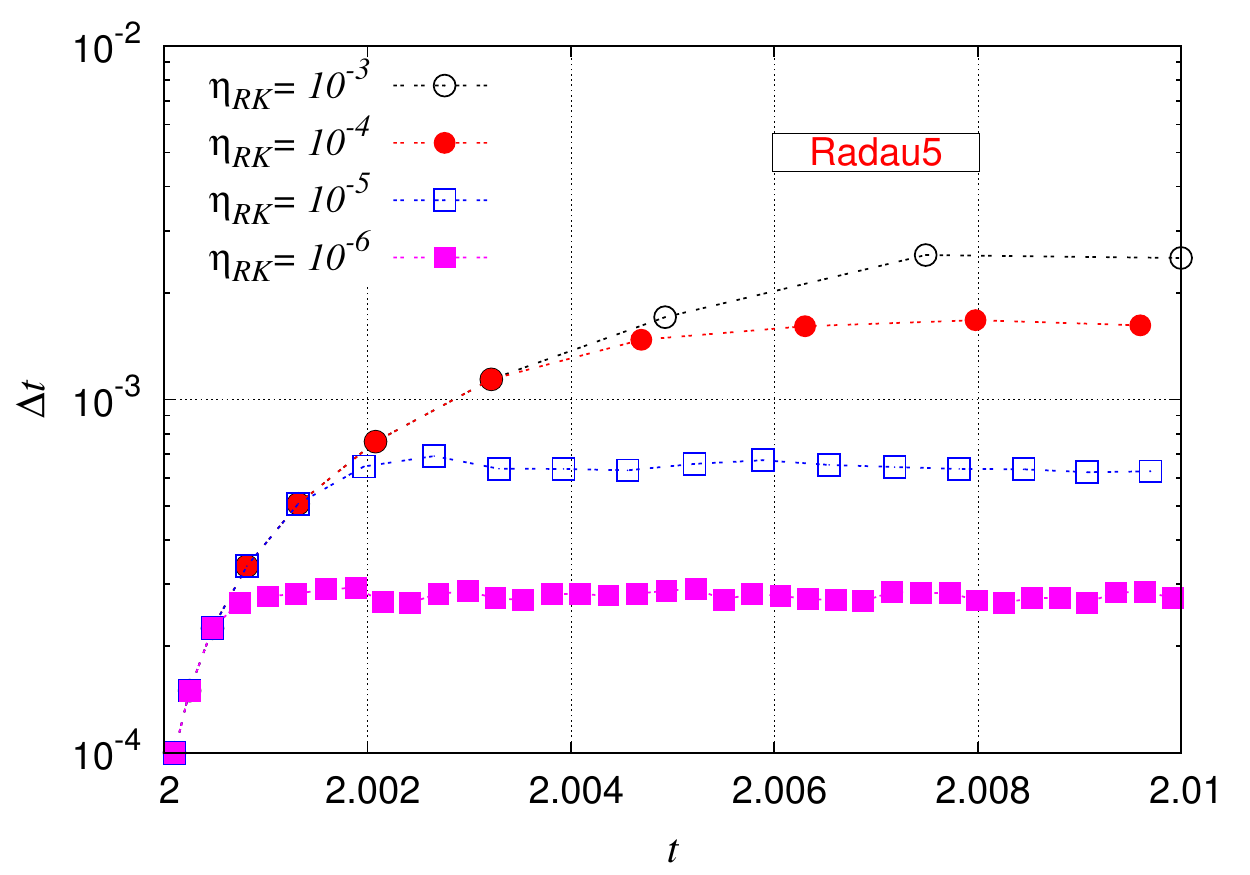}
\end{center}
\caption{Time--stepping with different accuracy tolerances, $\toltime$,
for the SDIRK4 (left) and Radau5 (right) schemes.}
\label{fig:BZ_2D_dt}
\end{figure}
\begin{table}[!htb]
\caption{Time integration with SDIRK4 for $t\in[2,2.01]$:
number of time steps, $n$;
maximum time step used, $\max \step_n$;
maximum number of Newton iterations, $\max k$;
maximum number of GMRES iterations, $\max \klin$;
CPU time in seconds.}
\label{TableSDIRK4}
\begin{center}
\begin{tabular}{l|c|c|c|c|c|c|c|} 
\cline{2-7}
\tlvs & 
\multirow{2}{*}{$\toltime$} & \multicolumn{5}{|c|}{SDIRK4} \\
\cline{3-7}
\tlvs & 
&$n$ & $\max \step_n$ & $\max k$ (per stage) & $\max \klin$ & CPU time (s) \\
\cline{2-7}
\tlvs & 
$10^{-3}$ &  $11$  & $1.95\times10^{-3}$ & $13$ ($3$)  & $15$  & $171.93$ \\
\tlvs & 
$10^{-4}$ &  $15$  & $1.01\times10^{-3}$ & $17$ ($4$)  & $13$  & $269.28$ \\
\tlvs & 
$10^{-5}$ &  $26$  & $4.65\times10^{-4}$ & $17$ ($4$)  & $11$  & $472.14$ \\
\tlvs & 
$10^{-6}$ &  $45$  & $2.51\times10^{-4}$ & $19$ ($4$)  & $10$  & $837.62$ \\
\cline{2-7}
\end{tabular}
\end{center}
\end{table}
\begin{table}[!htb]
\caption{Time integration with Radau5 for $t\in[2,2.01]$:
number of time steps, $n$;
maximum time step used, $\max \step_n$;
maximum number of Newton iterations, $\max k$;
maximum number of GMRES iterations, $\max \klin$;
CPU time in seconds.}
\label{TableRadau5}
\begin{center}
\begin{tabular}{l|c|c|c|c|c|c|} 
\cline{2-7}
\tlvs & 
\multirow{2}{*}{$\toltime$} & \multicolumn{5}{|c|}{Radau5} \\
\cline{3-7}
\tlvs & 
&$n$ & $\max \step_n$ & $\max k$& $\max \klin$ & CPU time (s) \\
\cline{2-7}
\tlvs & 
$10^{-3}$ &  $10$  & $2.56\times10^{-3}$  & $5$  & $60$  & $892.98$ \\
\tlvs & 
$10^{-4}$ &  $12$  & $1.64\times10^{-3}$  & $6$  & $56$  & $1268.11$ \\
\tlvs & 
$10^{-5}$ &  $19$  & $6.97\times10^{-4}$  & $6$  & $35$  & $2311.30$ \\
\tlvs & 
$10^{-6}$ &  $37$  & $2.98\times10^{-4}$  & $5$  & $21$  & $3416.22$ \\
\cline{2-7}
\end{tabular}
\end{center}
\end{table}

Figure~\ref{fig:BZ_2D_dt} shows the evolution of time steps
according to (\ref{eq:timestepping_dyn}), considering
$\alpha=1.5$ and $\step_0 = 10^{-4}$ at $t=2$ for
various accuracy tolerances: $\toltime$ between $10^{-3}$ and $10^{-6}$.
For this particular problem a roughly constant time step is
attained, consistent with the
quasi constant propagation speed of the wavefronts.
Tables~\ref{TableSDIRK4} and \ref{TableRadau5} gather information
on the performance of both solvers during
the time window $[2,2.01]$.
As also seen in Figure~\ref{fig:BZ_2D_dt}, larger time steps
for a given accuracy tolerance are used with the Radau5 scheme,
even though both schemes consider a third order, embedded method
to compute dynamically the integration time steps (\ref{eq2:time_stepping1_n}).
For SDIRK4, increasing the accuracy of the Newton solver involves more 
iterations even when smaller time steps are considered, showing
a rather low dependence on the time step size and a Newton solver
piloted mainly by its accuracy tolerance.
A different behavior is observed for Radau5 where smaller time steps
involve roughly the same number of Newton iterations, 
regardless of the Newton accuracy tolerance,
meaning that smaller time steps effectively improve the Newton solver.
All this is a direct consequence of the initialization of the Newton solver;
while for SDIRK4 the Newton solver at each stage is initialized using 
the previous stage solution and 
thus at some time within the current time step,
this is not the case for the present Radau5 solver
for which the larger the time step the worse the initial approximation.
In terms of the iterative linear solver, the number of iterations 
decreases considerably with smaller time steps even when tighter
convergence tolerances are considered. 
This is a direct consequence of the better
preconditioning of the more diagonal--dominant matrices in 
(\ref{eq2:Newton_simp}) and (\ref{eq2:Newton_simp_SDIRK}) for 
relatively small time steps.

In terms of CPU time, following
Tables~\ref{TableSDIRK4} and \ref{TableRadau5}
we can see that SDIRK4 is approximately $4$ to $5$ times faster than Radau5.
Updating the grid together with the multiresolution operations 
takes approximately $10$ to $13\,$\% for SDIRK4, whereas the time
integration, $84$--$90\,$\%. These numbers are within the range of 
values found in the literature for adaptive grid techniques
(see\eg \cite{Duarte_Phd}).
For Radau5 the multiresolution load goes
down to $2\,$\% with a roughly $98\,$\% of the CPU time 
allocated to the time integration, showing a clear
problem of performance.
There are two main reasons why a straightforward implementation
of Radau5 is not fully satisfactory.
First of all, a better initialization of the 
Newton solver is required to improve the convergence rate of the 
linear solver, regardless of the solver considered.
For example, in \cite{Hairer96} (\S~IV.8) all stages are initialized
by extrapolating from the previous time step and using an interpolation
polynomial based on the quadrature order conditions.
Even if an adaptive grid technique can considerably reduce the 
increase of data storage that the latter procedure involves for multi--dimensional PDEs,
we still have to introduce additional operations to initialize 
grid points that were not present during the previous time step.
The second problem is related to the size of the algebraic 
systems which are basically tripled in the case of Radau5.
The latter heavily impacts the performance of the linear
solver. In this particular implementation, most of the overload is
related to the preconditioning ILUT solver which was
implemented as a {\it black box}, contrary to the GMRES solver. 
A tailored ILUT solver implemented specifically for this data structure 
may have already improved its performance before considering parallel computing
implementations.

\subsection{Ignition model of diffusion flames}\label{subsec:Ignition}

We now consider the mathematical model
derived in \cite{Thevenin95} to investigate
the ignition dynamics of a diffusion flame, formed
while a reactive layer is being rolled--up
in a vortex.
The hydrodynamics is decoupled from species and energy transport equations by 
adopting a standard thermo--diffusive approximation, 
leading to a reaction--diffusion--convection
model.
A two--dimensional computational domain 
is considered
where 
pure and fresh hydrogen at 
temperature $T_{\F,0}$
initially occupies the 
upper half part, while
the remaining lower part of the domain is occupied
by hot air at $T_{\Ox,0}$.
By defining a Schvab--Zeldo'vich variable $Z$
and a reduced temperature $\theta$ given by
\begin{equation}\label{eq12:theta}
\theta = \frac{T-T_{\Ox,0}}{T_{\F,0}-T_{\Ox,0}},
\end{equation}
the mathematical model is given 
by a system of equations of the form \cite{Thevenin95}:
 \begin{equation}\label{eq12:eq_z_theta}
\left.
\begin{array}{l}
\ds
\partial _{t} Z + 
v_{x}  \partial _{x} Z +
v_{y}  \partial _{y} Z -
\left(\partial ^2_{x} Z  +
\partial ^2_{y} Z \right)= 0, \\[2.5ex]
\ds 
\partial _{t} \theta + 
v_{x}  \partial _{x} \theta +
v_{y}  \partial _{y} \theta -
\left(\partial ^2_{x} \theta  +
\partial ^2_{y} \theta \right)
= F(Z,\theta),
\end{array}
\right\}
\end{equation}
\begin{equation*}\label{eq12:F_Z_theta}
F(Z,\theta)=
{\rm Da}\, \phi \chi Y_{\Ox,0}
\left[
\frac{1-Z}{\phi \tau} + \frac{1}{\chi}(Z-\theta)
\right]
\left[
Z + \frac{\tau}{\chi}(Z-\theta)
\right]
 \e^{\left(- \tau_a/(1+\tau \theta)\right)},
\end{equation*}
with physical constant parameters:
${\rm Da}=1.65\times 10^{7}$,
$\phi=34.782608696$, 
$\chi=50$, 
$Y_{\Ox,0}=0.23$,
$\tau=-0.7$,
and $\tau_a=8$,
corresponding to $T_{\F,0}=300\,$K
and $T_{\Ox,0}=1000\,$K.
The velocity field 
$(v_{x},v_{y})$
is given by 
a single vortex centered on the planar interface between the two media,
which varies strongly in time and space.
Its tangential velocity is given by 
\begin{equation}\label{eq12:velocity_adim}
v_{\theta}(r,t)=
\ds \frac{{\rm Re}\, {\rm Sc}}{r}
\left(1-\e^{- r^2/(4\, {\rm Sc}\, t)} \right), 
\end{equation}
where $r(x,y)$ stands for the distance to the vortex center
$(x_0,y_0)=(0,0)$,
and with Reynolds and Schmidt numbers of
${\rm Re}=1000$ and ${\rm Sc}=1$, respectively.
In Cartesian coordinates, the velocity of a 
counter-clockwise rotating vortex
is thus given by
\begin{equation*}\label{eq12:velocity_cart}
v_{x} = \ds \left( \frac{y - y_{0} }{r}\right)
v_{\theta}, \quad
v_{y} = -  \ds \left( \frac{x - x_{0} }{r}\right)
v_{\theta}, \quad
r = \left[ (x - x_{0})^2 +  (y - y_{0})^2 \right]^{1/2}.
\end{equation*}

The physics of the phenomenon can be briefly described as follows.
A rotating vortex 
is introduced immediately at $t = 0$.
The resulting forced convection superposes to the diffusive mechanisms and
accelerates the mixture of the gases.
A diffusion flame then ignites along the contact surface of
both media, taking into account the important difference of temperatures
in those regions.
Once the flame is completely ignited, it 
propagates 
outwards from the center of the computational domain.
The complete phenomenon
encompasses thus very different physical regimes
like mixing, ignition, propagation,
which can be characterized
depending on the initial reactants configuration 
and on the imposed velocity field,
as studied in detail in \cite{Thevenin95}.

\subsubsection{High order temporal approximations}
\begin{figure}[!htb]
\begin{center}
\includegraphics[width=0.495\textwidth]{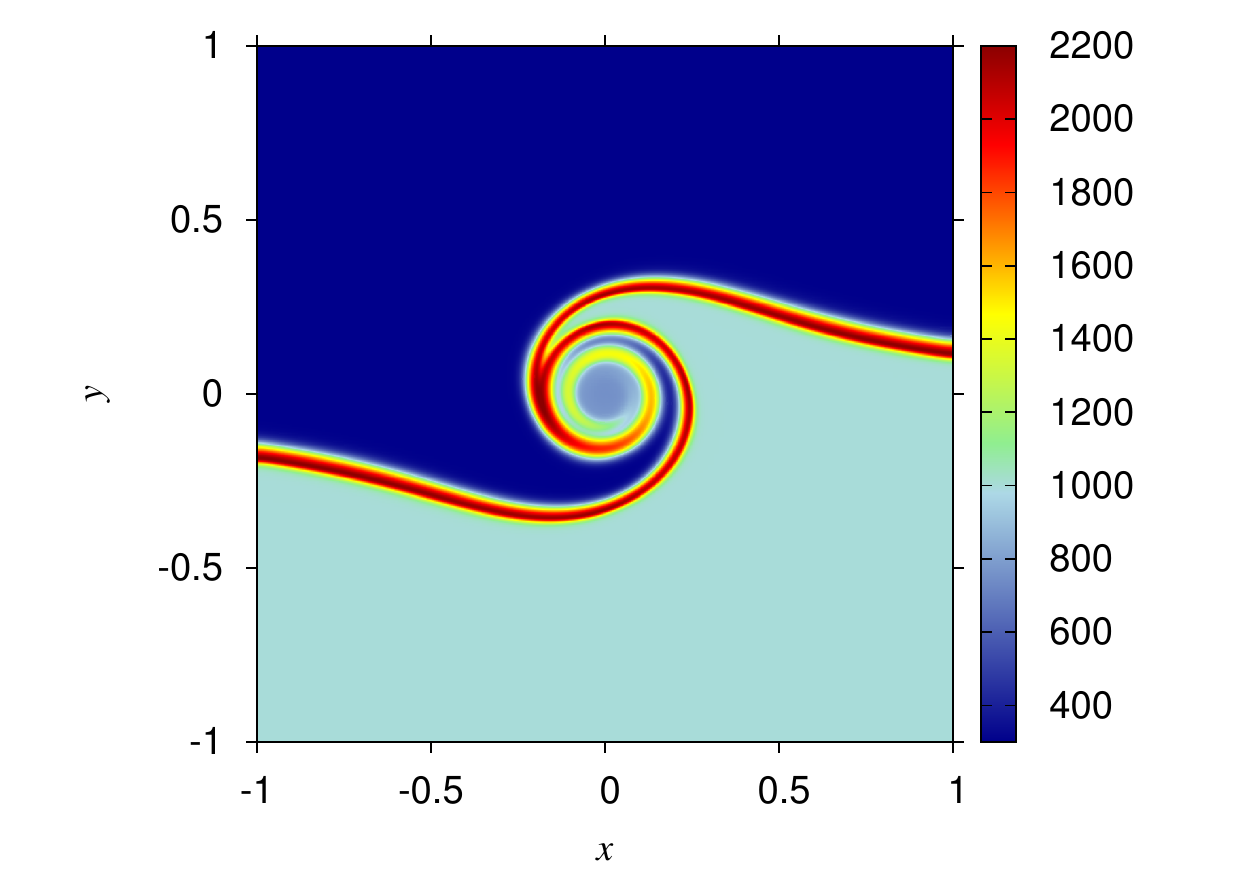}
\includegraphics[width=0.495\textwidth]{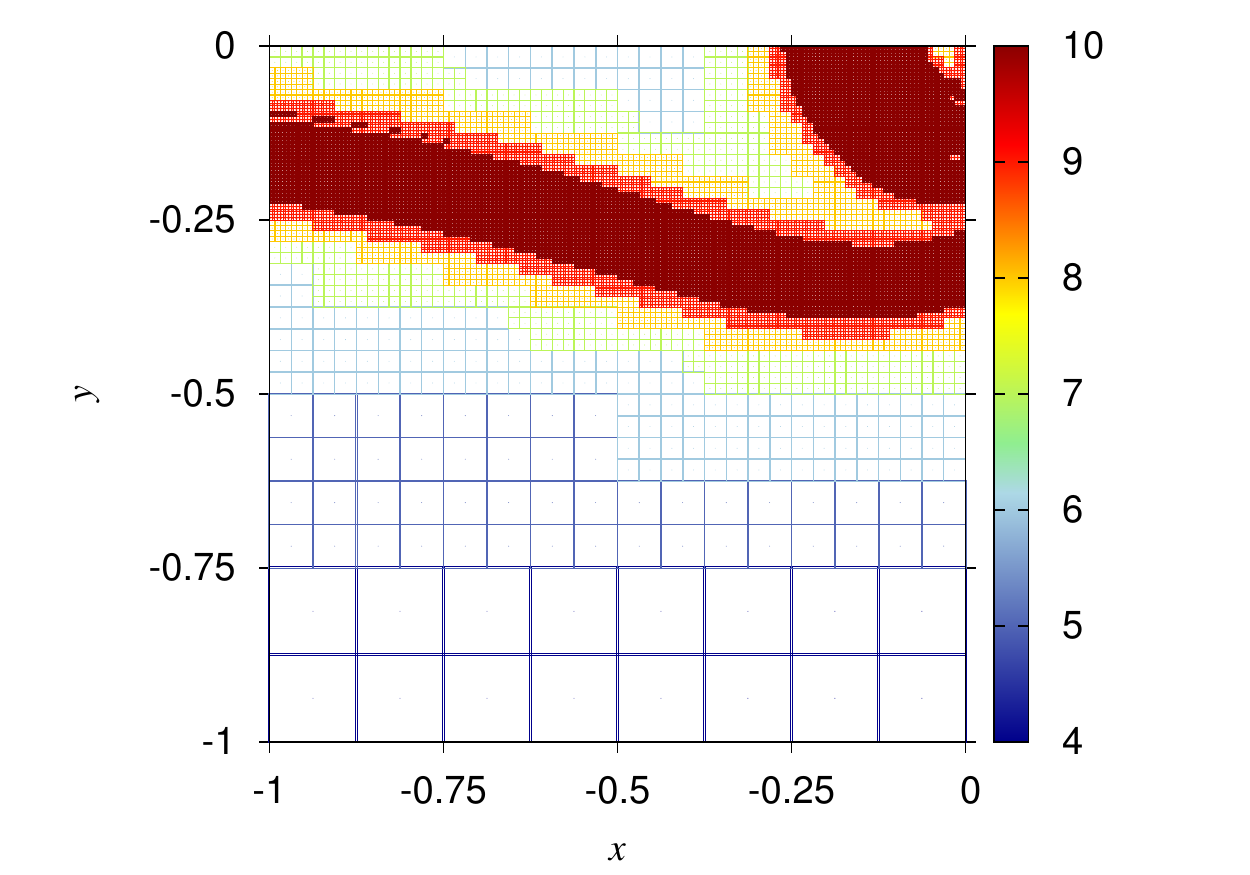}
\end{center}
\caption{Two--dimensional ignition model. Temperature $T$ deduced from (\ref{eq12:theta}) at $t=1.5\times 10^{-4}$ (left)
and the corresponding adapted grid for the zoomed region $[-1,0]\times[-1,0]$ (right).}
\label{fig:sol_Ignition_2D}
\end{figure}

We consider problem
(\ref{eq12:eq_z_theta}) in a two--dimensional configuration
with Neumann homogeneous boundary conditions,
using multiresolution analysis to adapt dynamically the 
spatial discretization grid.
The convective term is discretized in space
using a standard first--order upwind scheme.
As before the multiresolution analysis 
is parametrized as follows:
number of roots per direction,
$N_{{\rm R}x}=N_{{\rm R}y}=1$; 
maximum grid--level, $J=10$; and
accuracy tolerance,
$\tolMR = 10^{-3}$.
The finest grid has thus a spatial resolution of 
$1024 \times 1024$ over a computational domain of $[-1,1]\times [-1,1]$.
The model is simulated for a time window of $[0,1.5\times 10^{-4}]$,
using all the time integration solvers
previously described.
In all cases the following parameters were chosen:
$\knew = 30$, $\klinjac = \knew$, and
$\kappa = 10^{-2}$.
For this highly unsteady problem the number of active grid cells
increases from approximately 3\% of $1024^2$ for the initial
inert configuration, up to 13\% at the final time when the diffusion
flame is fully ignited along the contact surface
(see Figure~\ref{fig:sol_Ignition_2D}).
(See \cite{DuarteCFlame} for further details on the initialization of this problem.)
\begin{figure}[!htb]
\begin{center}
\includegraphics[width=0.49\textwidth]{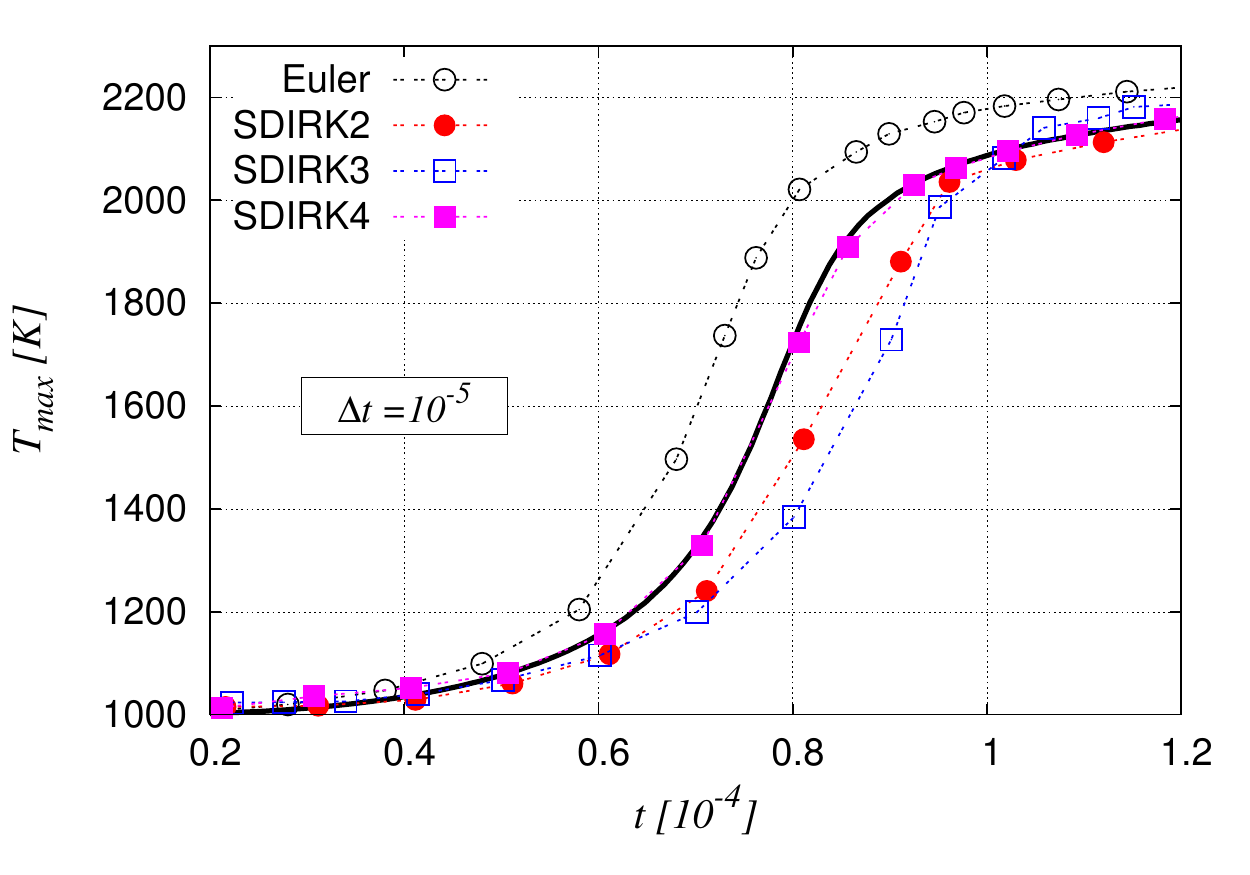}
\includegraphics[width=0.49\textwidth]{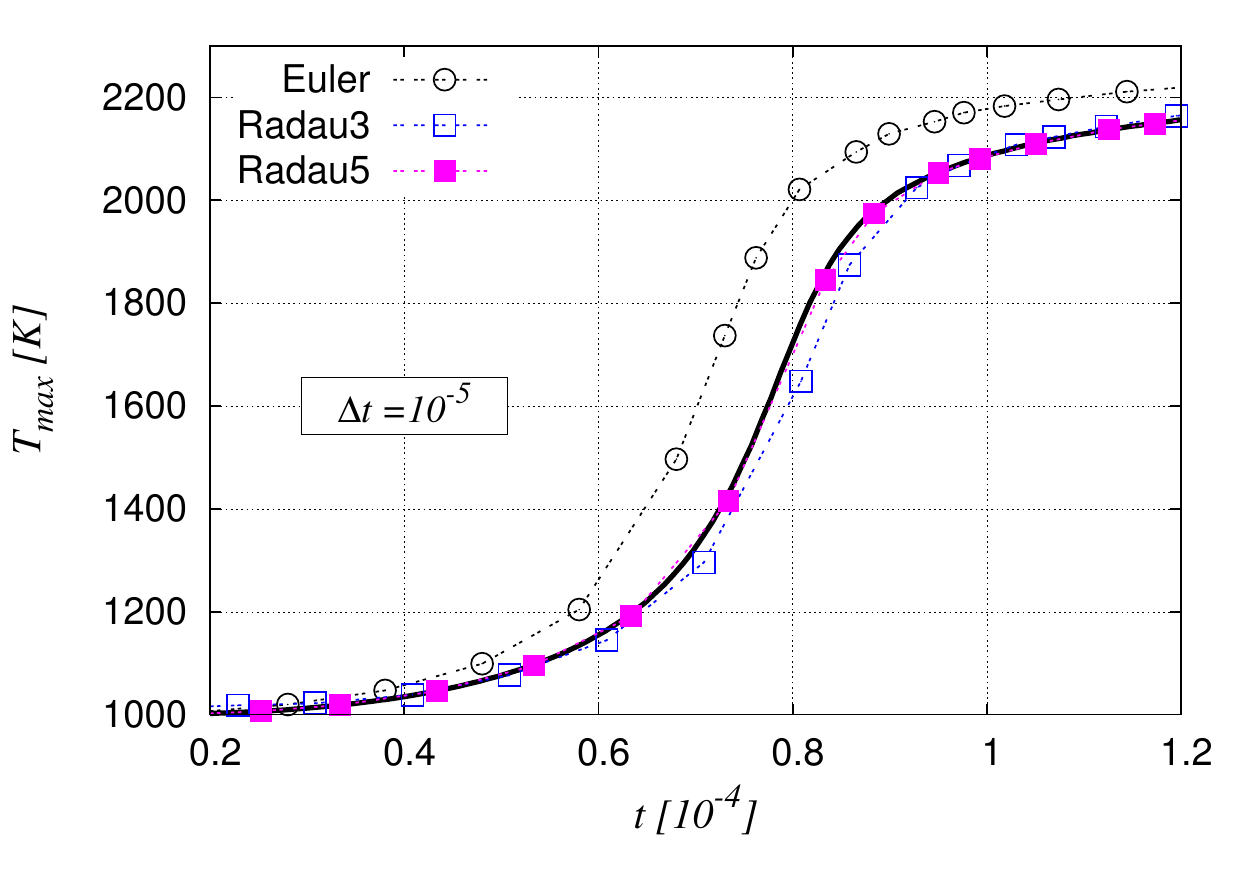}
\includegraphics[width=0.49\textwidth]{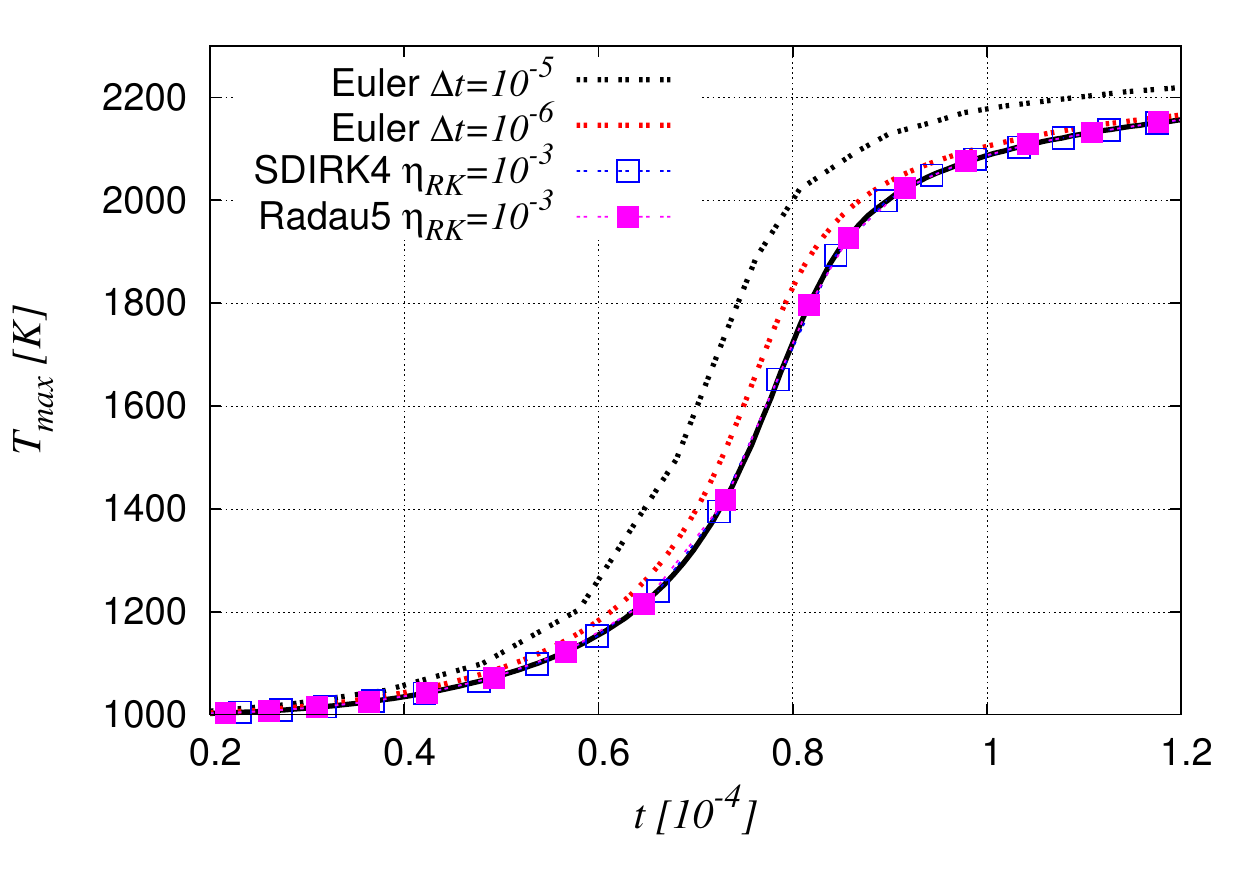}
\includegraphics[width=0.49\textwidth]{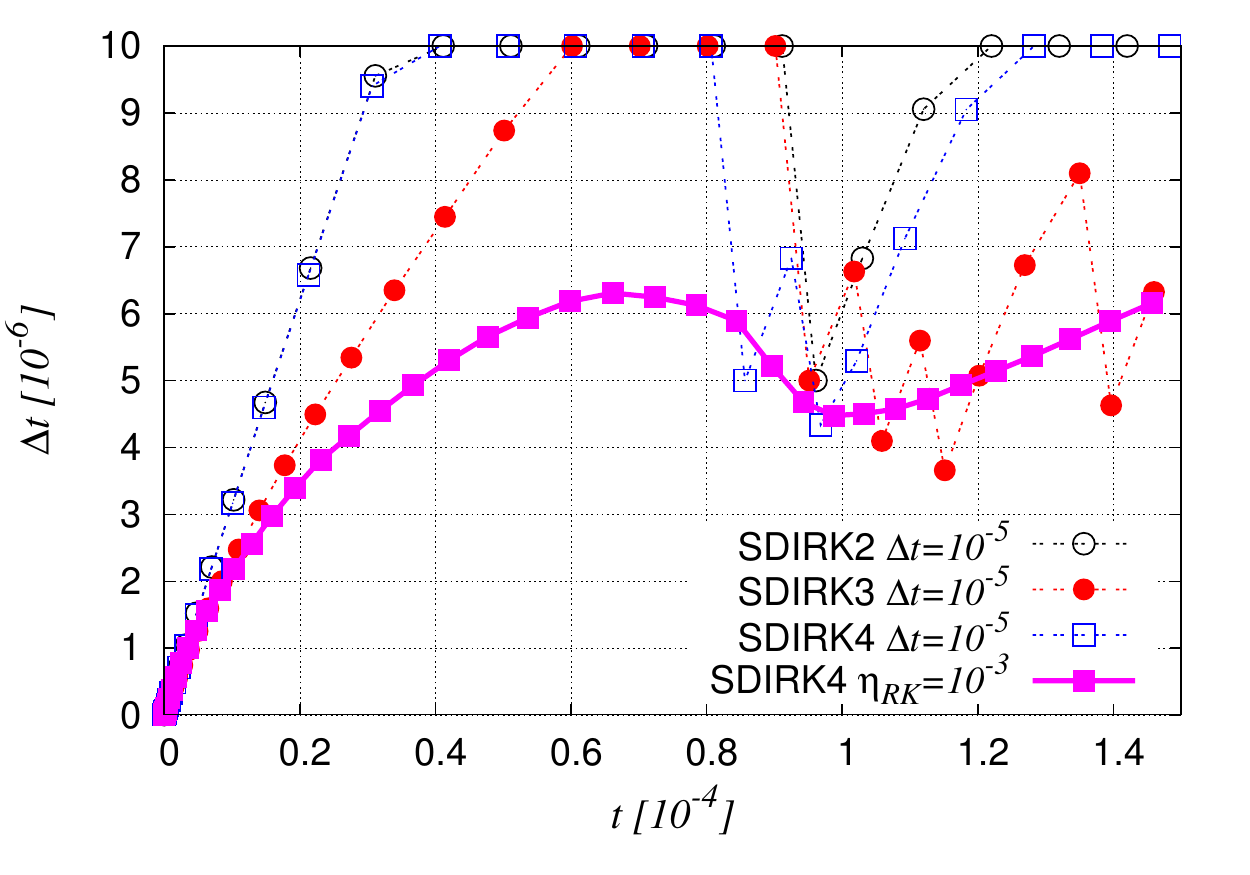}
\end{center}
\caption{Evolution of the maximum temperature $T_{max}$ using
a time step of $\step=10^{-5}$ for the SDIRK (top left) and 
Radau (top right) solvers. The solution computed with Radau5
and $\step=10^{-6}$ is depicted with a solid black line.
Similarly, time-adaptive solutions based on a tolerance of
$\toltime=10^{-3}$ are shown (bottom left).
The various time steps considered for the SDIRK solvers are 
also illustrated (bottom right); in all cases, $\step_0=10^{-8}$.}
\label{fig:res_Ignition}
\end{figure} 

Figure~\ref{fig:res_Ignition} shows the time evolution of the 
maximum temperature $T_{max}$ throughout the computational domain.
Notice that the initial $T_{max}$ corresponds to $T_{\Ox,0}=1000\,$K,
the hot air;
however, the fuel is initially at a much lower temperature
of $T_{\F,0}=300\,$K and thus the local temperature changes are in 
fact even more dramatic.
First, we consider a constant time step of $\step=10^{-5}$, but 
with the time--stepping procedure given by (\ref{eq:timestepping_const}). 
In all cases an initial time step of
$\step_0=10^{-8}$ was considered, taking into account that 
the velocity field (\ref{eq12:velocity_adim})
radically changes during the first time step.
A tolerance of $\tollin=10^{-5}$ is considered for the Newton solver.
Figure~\ref{fig:res_Ignition} (top) clearly shows the difference 
between the approximations obtained with different discretization
orders. As a reference solution we consider the one obtained with Radau5
and $\step=10^{-6}$.
The first order Euler solver introduces an ignition delay
of the order of one time step. This delay is subsequently corrected
by increasing the order of the time discretization with the same time step.
Considering the third order approximations, SDIRK3 and Radau3, the latter
performs much better for this particular problem because of its $L$--stability
capabilities,
given the strong transients that the stiff 
system (\ref{eq12:eq_z_theta}) models.
As a matter of fact, the $L$--stable SDIRK2 performs also better
than SDIRK3.
The difference of quality of the approximations can be further assessed
by considering the final temperature at $t=1.5\times 10^{-4}$,
once the ignition process is achieved and the strongest transients resolved.
Final $T_{max}$ goes from $2254.09\,$K for the Euler solver 
to $2227.96\,$K and $2224.74\,$K for SDIRK4 and Radau5, respectively.
For comparison, the reference Radau5 solution yields $2224.95\,$K.
With a smaller time step of $\step=10^{-6}$, the Euler solution still
shows a considerable difference 
(see Figure~\ref{fig:res_Ignition} (bottom left))
with a final $T_{max}$ of $2219.54\,$K.
All the other solvers with $\step=10^{-6}$ 
yield solutions with less than $0.3\,$K of difference
with respect to the Radau5 solution, except for SDIRK3 with about
$1\,$K of difference.

As an illustration, 
Figure~\ref{fig:res_Ignition} (bottom right)
also depicts the time steps considered for the various
SDIRK solvers.
Recalling that the time--stepping strategy is actually
influenced by the performance of the Newton and linear
solvers, we can see the impact
of strong physical changes during the ignition process
as all solvers need to use at some point smaller time steps.
Similar behaviors are observed for the Euler and Radau solvers
(not shown).
In particular it can be seen again how the non $L$--stable 
SDIRK3 is the most affected solver.
Notice that for these numerical experiments
we have chosen a relatively large $\knew$ 
as this allows for larger time steps 
according to (\ref{eq:def_nu});
consequently, $\klinjac$ is also large and the Jacobians
are never recomputed.
Time steps are thus reduced due to a bad convergence rate
of the Newton solver; a more conservative 
lower $\knew$ prevents this bad convergence rates 
since time steps would not even attain $10^{-5}$.
In general a careful tuning of parameters should be conducted with 
this ``constant'' time step strategy in order to 
get the best possible performance. This tuning can be 
highly problem--dependent.
That is why a time--stepping strategy based on an accuracy tolerance
is very convenient, as the
time steps can be effectively adapted to the various physical
scenarios within a prescribed accuracy, while reducing
the importance of the many parameters related to the 
Newton and linear solvers.
Time--adaptive solutions are shown in Figure~\ref{fig:res_Ignition} 
(bottom left) for SDIRK4 and Radau5 with $\toltime=10^{-3}$.
In terms of CPU time, with $\step=10^{-5}$
the Euler solver takes approximately $3.6$ minutes, 
compared to $5$ and $14.9$ for SDIRK4 and Radau5, respectively;
but the physics simulated with the first order method 
diverges considerably from the right one.
The time--adaptive SDIRK4 with $\toltime=10^{-3}$ takes 
approximately $4$ minutes, becoming a very promising alternative
to the cheaper but less accurate Euler scheme, especially
if one takes into account that a more accurate Euler solver with 
$\step=10^{-6}$ takes about $6$ minutes.

\section{Concluding remarks}\label{sec:conclusion}

We have considered high order, implicit integration schemes to solve stiff
multi--dimensional PDEs on adaptive multiresolution grids.
Such an adaptive technique yields highly compressed representations
within a user--prescribed accuracy tolerance, 
considerably reducing the computational requirements of 
implicit Runge--Kutta schemes.
In particular a competitive time--space adaptive strategy was introduced
to simulate models involving different physical scenarios with a
broad spectrum of time and space scales within a user--specified
level of accuracy.
By designing an appropriate procedure to evaluate functions
and represent linear systems 
within the multiresolution data structure, 
we have implemented several implicit Runge--Kutta schemes
of SDIRK-- and RadauIIA--type.
The resulting linear systems are completely independent of the grid 
generation or any other grid--related
data structure or geometric consideration.
Solving the algebraic problems constitute then a 
separate aspect from the multiresolution analysis itself,
while the same procedure remains perfectly valid for other
space adaptive techniques.

Three stiff models have been investigated to assess the 
computational performance of the numerical strategy in terms 
of accuracy and CPU time.
The computational analyses have thus proved 
that stiff PDEs can be effectively approximated 
with high order time discretization schemes 
with very limited computational resources.
In particular SDIRK schemes require roughly the same amount of memory
than a standard, low order Euler method.
More memory--demanding RadauIIA schemes can also be employed in
conjunction with adapted grids; however,
as previously discussed,
further enhancements 
are required to achieve better computational performances.
It was also shown that even in the presence of order reduction,
high order schemes yield more accurate solutions than low order ones.
The advantages of high order discretizations
have been especially highlighted when dealing with highly unsteady problems.
However, for problems of even larger size
parallel computing capabilities
must be developed within the current
context to achieve overall satisfactory results. 
Additionally, high order space discretization schemes,
well--suited for implicit schemes \cite{Noskov2005,Noskov2007,Dobbins2010},
could be also considered in conjunction with grid adaptation 
to further enhance the computational performance. 
These issues constitute particular topics of our
current research.


\appendix

\section{Details on multiresolution analysis}\label{app:multiresolution}
Defining $\Omega_{\lambda}:=\Omega_{j,k}$, 
we denote $|\lambda|:=j$ if $\lambda \in S_j$, while subscript $k\in \Delta_j \subset \Z^d$
corresponds to the position of the cell within $S_j$. 
For instance, 
in Cartesian coordinates
the univariate dyadic intervals in $\R$ are given by
\begin{equation}\label{eq3:dyadic_1D}
\Omega_{\lambda}=
\Omega_{j,k} := ]2^{-j}k,2^{-j}(k+1)[,\ \lambda \in S_j:= \{(j,k) \ \mathrm{s.t.} \ j\in (0,1,\ldots,J), \, k\in \Z \},
\end{equation}
and the same follows for higher dimensions.

Following \cite{Cohen03},
the projection operator $P^j_{j-1}$
maps $\disc U_j$ to $\disc U_{j-1}$.
It is obtained
through exact a\-ve\-ra\-ges computed at the finer level by
\begin{equation}\label{eq3:projection}
  u_{\lambda} = |\Omega_{\lambda}|^{-1} 
\sum_{|\mu|=|\lambda| +1,\Omega_{\mu} 
\subset \Omega_{\lambda}} |\Omega_{\mu}|u_{\mu}.
 \end{equation}
As far as grids are nested, this
projection operator is {\it exact} and {\it unique}
\cite{cohen2000a}.
On the other hand,
the prediction operator $P^{j-1}_j$
maps $\disc U_{j-1}$ to 
an approximation $\widehat{\disc U}_j$ of $\disc U_{j}$.
There are several choices to define $P^{j-1}_j$, but 
two basic constraints are usually imposed:
\begin{enumerate}
\item The prediction is local\ie $\widehat{u}_{\lambda}$ depends on the values $u_{\mu}$
 on a finite stencil $R_{\lambda}$ surrounding $\Omega_{\lambda}$, where $|\lambda|=|\mu| +1$.
\item The prediction is {\it consistent} with the projection in the sense that
 \begin{equation}\label{eq3:prediction_consist}
   u_{\lambda} = |\Omega_{\lambda}|^{-1}\sum_{|\mu|=|\lambda| +1,\Omega_{\mu} \subset \Omega_{\lambda}} |\Omega_{\mu}|\widehat{u}_{\mu};
 \end{equation}
    \textit{i.e.}, one can retrieve the coarse cell averages from the
predicted values: 
\begin{equation}
 P_{j-1}^j \circ P_j^{j-1} = \I.
\end{equation}
\end{enumerate}

With these operators, for each cell $\Omega _{\lambda}$
the prediction error or {\it detail} is defined
as the difference between the exact and predicted values,
\begin{equation}\label{eq3:MR_detail}
d_{\lambda} := u_{\lambda} - \widehat{u}_{\lambda},
\end{equation}
or in terms of inter--level operations:
$d_{\lambda} = u_{\lambda} - P_{|\lambda|}^{|\lambda|-1} \circ P_{|\lambda|-1}^{|\lambda|} u_{\lambda}$.
The consistency property (\ref{eq3:prediction_consist}) and the 
definitions of the
projection operator (\ref{eq3:projection}) and 
of the detail (\ref{eq3:MR_detail}) imply
\begin{equation}\label{eq3:MR_consis_2}
\sum_{|\lambda|=|\mu| +1,\Omega_{\lambda} \subset \Omega_{\mu}} |\Omega_{\lambda}|d_{\lambda}=0.
 \end{equation}
We can then construct as shown in \cite{Cohen03}
a {\it detail vector} defined as $\disc D_j:=(d_{\lambda})_{\lambda \in \nabla_j}$,
where the set $\nabla_j \subset S_j$ is obtained by removing for
each $\mu \in S_{j-1}$ one $\lambda \in S_j$ ($\Omega_{\lambda} \subset \Omega_{\mu}$)
in order to avoid redundancy (considering (\ref{eq3:MR_consis_2}))
and to get the one--to--one correspondence (\ref{eq3:one_one_cor}).

Given a set of indices $\Lambda \subset \nabla^J$,
where $\nabla^J := \bigcup_{j=0}^J \nabla_j$ with 
$\nabla_0 := S_0$,
the thresholding operator $\thr_{\Lambda}$
is such that leaves unchanged the components $d_{\lambda}$ 
of the multi--scale representation $\disc M_J$ in (\ref{eq3:M_U_M})
if $\lambda \in \Lambda$,
and replaces it by $0$ otherwise.
Defining the level--dependent threshold values 
$(\epsilon_0,\epsilon_1,\ldots,\epsilon_J)$,
the set $\Lambda$ is given by 
\begin{equation}\label{eq3:MR_Lambda}
\lambda \in \Lambda \ {\rm if} \ 
\|d_{\lambda}\|_{L^p} \geq \epsilon_{|\lambda|}.
\end{equation}
Applying $\thr_{\Lambda}$ on the multi--scale decomposition 
$\disc M_J$ of $\disc U_J$ 
amounts then to building the multiresolution
approximation $\disc U_J^{\epsilon} := \adap_{\Lambda}\disc U_J$ to $\disc U_J$, 
where the operator 
$\adap_{\Lambda}$ is given by
\begin{equation*}\label{eq3:adap_operator_def}
\adap_{\Lambda}:=\MR^{-1}\thr_{\Lambda}\MR,
\end{equation*}
in which all details below a certain level of regularity
have been discarded.
The bound (\ref{eq3:adap_error_eps}) is thus verified 
with the level--dependent threshold values:
\begin{equation*}\label{eq3:epsilon_j}
\epsilon_j = 2^{d(j-J)/2}\tolMR, \quad j=0,1,\ldots,J.
\end{equation*}

\section{Details on implicit Runge--Kutta schemes}\label{app:IRK}
The solution $\disc U(t_0+\step)$ of problem (\ref{eq:gen_disc_prob})
is approximated by $\disc U_1$, computed as
\begin{align}
&\disc g_i = \disc U_0 + \step \ds \sum_{j=1}^{s} a_{ij} \disc F \left(t_0+c_j \step, \disc g_j\right),
\qquad i=1,\ldots,s; \label{eq2:runge_kutta_1}\\
&\disc U_1= \disc U_0+ \step \ds \sum_{j=1}^s b_j \disc F \left(t_0+c_j \step, \disc g_j\right), 
\label{eq2:runge_kutta_2}
\end{align}
where the time dependence for $\disc F(\disc U (t))$ was added for the sake of clarity. 
Following the approach established in \cite{Hairer96},
the set of arrays
$\disc z_1,\ldots,\disc z_s$,
are defined such that 
$\disc z_i =\disc g_i-\disc U_0$,
$i=1,\ldots,s$,
and hence,
\begin{equation}\label{eq2:runge_kutta_radau}
\disc z_i =\step \ds \sum_{j=1}^s a_{ij} \disc F(t_0+c_j\step,\disc U_0+ \disc z_j),
\qquad i=1,\ldots,s.
\end{equation}
Therefore, knowing the solution $\disc z_1,\ldots,\disc z_s$ 
implies an explicit formula for $\disc U_1$ in (\ref{eq2:runge_kutta_2}), for which 
$s$ additional function evaluations are required.
These extra computation can nevertheless be avoided
if the matrix $\vec A$ is nonsingular, as seen in (\ref{eq2:u1_radau}).
 
A standard iterative Newton solver for system
(\ref{eq2:nonlinear_sys}) 
amounts to solving at each iteration a linear system 
that requires the inversion of the block--matrix:
\begin{equation}\label{eq2:Newton_linearsys}
\left(
\begin{array}{ccc}
\Ibf_{m \times N} -\step a_{11}  
\disc J_1
(\disc z_1)
&\ldots
&-\step a_{1s} 
\disc J_s
(\disc z_s)
\\
\vdots&\ddots&\vdots\\
-\step a_{s1} 
\disc J_1
(\disc z_1)
&\ldots
&\Ibf_{m \times N} - \step a_{ss}  
\disc J_s
(\disc z_s)
\end{array}\right),
\end{equation}
where $\disc J_i(\disc z_i) := \partial_{\disc U} \disc F(t_0+c_i \step,\disc U_0+ \disc z_i)$,
$i=1,\ldots,s$,
stands for the Jacobian $\partial_{\disc U} \disc F := \partial \disc F(\disc U)/\partial \disc U$
of size $m \times N$, evaluated at the various inner stages.
Following \cite{Hairer96},
the {\it simplified} Newton solver (\ref{eq2:Newton_simp})
approximates all Jacobians in (\ref{eq2:Newton_linearsys}) as
\begin{equation}\label{eq2:Jac_approx}
\disc J_i(\disc z_i) \approx \disc J_0 := \partial_{\disc U} \disc F (t_0,\disc U_0),
\quad i=1,\ldots,s.
\end{equation}

\section{Butcher tableau of IRK schemes}\label{app:Butcher_tab}
We consider a two--stage SDIRK scheme \cite{Hairer87} (Table II.7.2) given by
\begin{equation}\label{eq:SDIRK23}
\begin{tabular}{c|c c }
$\gamma$ & $\gamma$ &   \\
$1-\gamma$ & $1-2\gamma$ & $\gamma$  \\
\hline
& $1/2$ & $1/2$ 
\end{tabular}
\qquad
d_1 = \frac{3\gamma - 1}{2\gamma^2}, 
\quad
d_2 = \frac{1}{2\gamma}.
\end{equation}
With $\gamma = (2 \pm \sqrt{2})/2$
the method (\ref{eq:SDIRK23})
is of second order 
($p=2$, $s=p$) and $L$--stable.
For $\gamma = (3 + \sqrt{3})/6$
it becomes a third order scheme 
($p=3$, $s=p-1$), 
but is only $A$--stable.
The fourth order SDIRK4 scheme proposed
in \cite{Hairer96} (Table IV.6.5) with $s=5$, is
\begin{equation}\label{eq:SDIRK4}
\begin{tabular}{c|c c c c c}
$\ds \frac{1}{4}$ & $\ds \frac{1}{4}$ & & & & \\[1.5ex]
$\ds \frac{3}{4}$ & $\ds \frac{1}{2}$ & $\ds \frac{1}{4}$ & &  &  \\[1.5ex]
$\ds \frac{11}{20}$ & $\ds \frac{17}{50}$ & $- \ds \frac{1}{25}$ & $\ds \frac{1}{4}$ &  &  \\[1.5ex]
$\ds \frac{1}{2}$ & $\ds \frac{371}{1360}$ & $- \ds \frac{137}{2720}$ & $\ds \frac{15}{544}$ & $\ds \frac{1}{4}$ & \\[1.5ex]
$1$ & $\ds \frac{25}{24}$ & $- \ds \frac{49}{48}$ & $\ds \frac{125}{16}$ & $-\ds \frac{85}{12}$ & $\ds \frac{1}{4}$  \\[1.25ex]
\hline\\[-2.ex]
& $\ds \frac{25}{24}$ & $- \ds \frac{49}{48}$ & $\ds \frac{125}{16}$ & $-\ds \frac{85}{12}$ & $\ds \frac{1}{4}$ 
\end{tabular}
\end{equation}
The latter is an $L$--stable and stiffly accurate method.

The RadauIIA methods of order 3 and 5 are given,
respectively, by
(\cite{Hairer96} Tables IV.5.5 and IV.5.6)
\begin{equation}
\begin{tabular}{c|cc}
$\ds \frac{1}{3}$&$\ds\frac{5}{12}$&$- \ds\frac{1}{12}$\\[1.5ex]
$1$&$\ds\frac{3}{4}$&$\ds\frac{1}{4}$\\[1.25ex]
\hline\\[-2ex]
&$\ds\frac{3}{4}$&$\ds\frac{1}{4}$
\end{tabular}
\end{equation}
and
\begin{equation}\label{eq:Radau5}
\begin{tabular}{c|ccc}
$\ds \frac{4-\sqrt{6}}{10}$&$\ds\frac{88-7\sqrt{6}}{360}$&$\ds\frac{296-169\sqrt{6}}{1800}$&$\ds\frac{-2+3\sqrt{6}}{225}$\\[1.5ex]
$\ds \frac{4+\sqrt{6}}{10}$&$\ds\frac{296+169\sqrt{6}}{1800}$&$\ds\frac{88+7\sqrt{6}}{360}$&$\ds\frac{-2-3\sqrt{6}}{225}$\\[1.5ex]
$1$&$\ds\frac{16-\sqrt{6}}{36}$&$\ds\frac{16+\sqrt{6}}{36}$&$\ds\frac{1}{9}$\\[1.25ex]
\hline\\[-2ex]
&$\ds\frac{16-\sqrt{6}}{36}$&$\ds\frac{16+\sqrt{6}}{36}$&$\ds\frac{1}{9}$
\end{tabular}
\end{equation}

\section{Numerical computation of Jacobians}\label{app:Jac}
The Jacobian 
at a given stage $s$:
$\disc J_s = \partial_{\disc U} \disc F (\disc U_s)= (J^s_{i,j})_{i,j\in {\rm I}^n_\leaf}$,
needs to be numerically approximated at each leaf of the adapted grid.
The latter is done by considering for each $i=1, 2, \ldots, N^n_\leaf$, $\lambda = h_n^{-1}(i)$, 
the following expression
\begin{equation}
J^s_{h_n(\lambda),h_n(\mu)} = 
\frac{F_{\lambda}(\disc U_s+ \delta_{\mu} \One_{\mu} )-F_{\lambda}(\disc U_s)}
{\delta_{\mu}},
\qquad \forall
\mu \in 
R_{\varPhi}[\lambda],
\end{equation}
with the one--dimensional array $ \One_{\mu} \in R^{N^n_\leaf}$,
such that $\mathbf{1}_{\mu}(i)=1$ 
for $i=h_n(\mu)$ and
$\mathbf{1}_{\mu}(i)=0$, otherwise;
$\delta_{\mu}$ is a small perturbation
taken here as
$\delta_{\mu}=(10^{-16}\times \max(10^{-5},|u^s_{\mu}|))^{0.5}$,
following \cite{Hairer96}.
Recalling (\ref{eq:disc_prob_local_flux})
and the conservation property, the latter 
perturbation is limited in practice to the stencil 
$R_{\varPhi}^+[\lambda]$. 
That is, 
perturbed fluxes $\varPhi^+_{\lambda}(\disc U_s+\delta_{\mu}\One_{\mu})$
are evaluated for all $\mu \in R_{\varPhi}^+[\lambda]$ 
to compute $J^s_{h_n(\lambda),h_n(\mu)}$,
while  
$-\varPhi^+_{\lambda}(\disc U_s+\delta_{\mu}\One_{\mu})$
is used to compute the Jacobian entry at cell $\mu$ due to perturbed cell $\lambda$:
$J^s_{h_n(\mu),h_n(\lambda)}$.
Notice that only cells contained in the adapted grid\ie
the leaves, are perturbed to compute the Jacobians.
Close to level interfaces where a ghost--cell
$\Omega_{\nu}$
might be
contained in the flux stencil,
we perturb the latter with the same $\delta_{\mu}$
corresponding to its parent 
$\Omega_{\mu} \supset \Omega_{\nu}$,
which is of course a leaf.
Geometric proportions must then be taken into account when
transferring fluxes between cells of different resolution.

\section{Embedded Runge--Kutta schemes}\label{app:embeddedRK}
The error measures used 
in this work rely upon
the lower order embedded IRK schemes
introduced in \cite{Hairer96} for the SDIRK4 and Radau5 schemes.
For SDIRK4, the error approximation $err$ is computed from
\begin{equation}
\hat{\disc U}_1 - \disc U_1 = \sum_{i=1}^5 e_i \disc z_i,
\end{equation}
with coefficients 
$\vec e^T := (e_1,\ldots,e_5) = (23/6,17/12,-125/4,85/3,1)$,
which is accurate to third order ($\hat{p}=3$).
In the case of Radau5, we consider
\begin{equation}
\hat{\disc U}_1 - \disc U_1 = K \step \disc F (t_0,\disc U_0)+
\sum_{i=1}^3 e_i \disc z_i,
\end{equation}
with $K=10^{-1}$,
$\vec e^T = K (-13-7\sqrt{6},-13+7\sqrt{6},-1)$,
which is also accurate to third order ($\hat{p}=3$).
As pointed out in 
\cite{Kvaerno2004}
neither of these error estimates has $A$--stability properties.
A remedy was proposed in \cite{Hairer96} for their Radau5 solver,
that is currently not implemented in our code.

%
\bibliographystyle{plain}
\bibliography{biblio_Implicit}

\end{document}